\documentclass[12pt,reqno]{amsart}

\usepackage{amsmath}
\usepackage{amssymb}
\usepackage{amsthm}
\usepackage[enableskew]{youngtab}
\usepackage{qtree}
\usepackage{eepic}
\usepackage{epsfig,graphicx}
\usepackage{color}

\allowdisplaybreaks

\newcommand{\red}[1]{\color{red}{#1}\color{black}}

\newtheorem{theorem}{Theorem}[section]
\newtheorem{prop}[theorem]{Proposition}
\newtheorem{lemma}[theorem]{Lemma}
\newtheorem{cor}[theorem]{Corollary}

\theoremstyle{remark}
\newtheorem*{example}{Example}
\newtheorem*{remark}{Remark}
\newtheorem*{note}{Note}

\theoremstyle{definition}
\newtheorem*{definition}{Definition}
\newtheorem*{claim}{\sc Claim}

\def\l{\lambda}
\def\m{\mu}
\def\n{\nu}
\def\a{\alpha}
\def\b{\beta}

\begin{document}

\newbox\Adr
\setbox\Adr\vbox{
\centerline{\sc Cristina M. Ballantine$^{1*}$ and Rosa C. Orellana$^{2**}$}
\vskip18pt
\centerline{$^1$College of the Holy Cross, Worcester, MA 01610, USA}
\centerline{Email: {\tt cballant@holycross.edu}}
\vskip6pt
\centerline{$^2$Dartmouth College, Hanover, NH 03755, USA}
\centerline{Email: {\tt Rosa.C.Orellana@Dartmouth.edu}}
}

\title{A Combinatorial Interpretation for the Coefficients
in the Kronecker Product $s_{(n-p,p)}*s_{\lambda}$}
\author[Cristina M. Ballantine and Rosa C. Orellana]{\box\Adr}

\thanks{\hphantom{$^*$}$^*$Partially supported by the Fulbright Commission}
\thanks{$^{**}$Partially supported by the Wilson Foundation}

\begin{abstract}
In this paper we give a combinatorial interpretation for the coefficient of
$s_{\nu}$ in the Kronecker product $s_{(n-p,p)}\ast s_{\l}$, where $\l=(\l_1,
\ldots, \l_{\ell(\l)})\vdash n$,  if   $\ell(\l)\geq 2p-1$ or $\l_1\geq 2p-1$;
that is, if $\l$ is not a partition inside the $2(p-1)\times 2(p-1)$ square.
For $\lambda$ inside the square our combinatorial interpretation provides
an upper bound for the coefficients. In general, we are able to combinatorially
compute these coefficients for all $\l$ when $n>(2p-2)^2$. We use this combinatorial
interpretation to give characterizations for multiplicity free Kronecker products.  We have
also obtained some formulas for special cases.
\end{abstract}

\maketitle

%
%

\section*{Introduction}
\pagenumbering{arabic}
\addtocounter{page}{0}
\markboth{}{}

Let $\chi^{\l}$ and $\chi^{\m}$ be the irreducible characters of $S_n$ (the
symmetric group on $n$ letters) indexed by the partitions $\l$ and $\m$ of $n$. The
\emph{Kronecker product} $\chi^{\l}\chi^{\m}$ is defined by
$(\chi^{\l}\chi^{\m})(w)=\chi^{\l}(w)\chi^{\m}(w)$ for $w\in S_n$. Then
$\chi^{\l}\chi^{\m}$ is the character that corresponds to the diagonal action of
$S_n$ on the tensor product of the irreducible representations indexed by $\l$ and
$\m$. We have
$$\chi^{\l}\chi^{\m} =\sum_{\nu\vdash n}g_{\l,\m,\n}\chi^{\nu},$$
where $g_{\l,\m,\n}$ is the multiplicity of $\chi^{\nu}$ in $\chi^{\l}\chi^{\m}$.
Hence the  coefficients $g_{\l,\m,\n}$ are non-negative integers. 

By
means of the Frobenius map one defines the Kronecker (internal) product on the Schur
symmetric functions by
$$s_{\l}\ast s_{\m}=\sum_{\n\vdash n} g_{\l,\m,\n} s_{\n}.$$
A formula for decomposing the Kronecker product is unavailable, although the problem
has been studied for nearly one hundred years. In recent years Lascoux \cite{la},
Remmel \cite{r1,r2}, Remmel and Whitehead \cite{rwh} and Rosas \cite{ro} derived
closed formulas for  Kronecker products of Schur functions indexed by two row shapes
or hook shapes. Gessel \cite{ge} obtained a combinatorial interpretation for zigzag
partitions. However, a combinatorial interpretation is still lacking even in the
case when both $\lambda$ and $\mu$ are two row partitions.

The objective
of this paper is to provide a combinatorial interpretation for the Kronecker
coefficients comparable to the Littlewood-Richardson rule which is defined in terms
of the so-called Littlewood-Richardson tableaux (see Section~1 for the definition).
In this paper we give a combinatorial interpretation for the coefficient of $s_{\n}$
in $s_{(n-p,p)}\ast s_{\l}$, if $\l_1\geq 2p-1$ or $\ell(\l) \geq 2p-1$, in terms of
what we call \emph{Kronecker Tableaux}. In particular, our combinatorial
interpretation holds for all $\l$ if $n>(2p-2)^2$.  For a  general $\l$, the number
of Kronecker tableaux always gives an upper bound for the Kronecker coefficients.
Furthermore, using our combinatorial rule we obtain that $g_{(n-p,p),\l,\n}=0$
whenever the intersection of $\l$ and $\n$ has less than $p$ boxes. 

The
techniques we use to obtain our main theorem are purely combinatorial and rely both
on the Jacobi-Trudi identity and on the Garsia--Remmel rule \cite{gr} for decomposing
the Kronecker product of a homogeneous symmetric function and a Schur symmetric
function. 


%
%

One can easily deduce from existing formulas that
$s_{(n-1,1)}\ast s_\l$ is multiplicity free if and only if $\l=(a^k, b^l)$ or
$\l=(a^k)$, where $a,k,b,l$ are non-negative integers. If $p \geq 2$, we have used
our combinatorial rule to determine the partitions $\l$ for which the Kronecker
product $s_{(n-p,p)}\ast s_\l$ is multiplicity free. We have determined that if
$n\geq 6$, $s_{(n-2,2)}\ast s_{\l}$ is multiplicity free if and only if $\l =(n),
(1^n), (n-1,1), (2,1^{n-2})$ or $\l=(a^k)$, where $a,k$ are non-negative integers.
If $n \geq 16$,  $s_{(n-3,3)}\ast s_\l$ is multiplicity free if and only if $\l
=(n), (1^n), (n-1,1), (2,1^{n-2})$ and if $n$ is even also $\l=(n/2,n/2)$ or
$\l=(2^{n/2})$. If $p\geq 4$ and $n > (2p-2)^2$, then $s_{(n-p,p)} \ast s_\l$ is
multiplicity free if and only if $\l =(n), (1^n), (n-1,1),
(2,1^{n-2})$. 

Other applications of our combinatorial interpretation for $g_{(n-p,p),\l,\n}$
include formulas for some special partitions $\l$ and $\n$.  The formulas obtained
do not have cancellations and are easy to program.

The paper is organized as follows. In Section~1 we give preliminary
definitions and set the notation used throughout the paper. In Section~2 we give a
variation of the Remmel--Whitney \cite{rwy} algorithm for expanding the skew Schur
function $s_{\l/\m}$. We then use this algorithm to prove that the symmetric
function $s_{\l/\a}s_{\a} -s_{\l/\b}s_{\b}$, where $\b=(\a_1-1,\a_2,\ldots,
\a_{\ell(\a)})$, is Schur positive if and only if $\l_1\geq 2\a_1-1$.  In Section~3
we define the Kronecker tableaux and give the combinatorial interpretation for
$g_{(n-p,p),\l,\n}$. In the last section we apply our combinatorial rule to give
characterizations for multiplicity free  $s_{(n-p,p)} \ast s_{\l}$. We  also give
closed formulas for several special cases.  For instance, we
 give a general formula for the coefficient of $s_{(n-t,t)}$
 in the product $s_{(n-p,p)}\ast s_{\lambda}$ and
show that these coefficients are unimodal for some special cases of $\lambda$.
\medskip

\noindent
{\bf Acknowledgement.} The authors are grateful to Christine Bessenrodt for useful
suggestions and comments.

\section{Preliminaries and Notation}
Details and proofs for the contents of this section can be found in \cite[Chap.
7]{s}. A \emph{partition} of a non-negative integer $n$ is a weakly decreasing
sequence of non-negative integers, $\l:=(\l_1,\l_2,\cdots,\l_{\ell})$, such that
$|\l|=\sum \l_i=n$. We write $\l\vdash n$ to mean $\l$ is a partition of $n$. The
nonzero integers $\lambda_i$ are called the \emph{parts} of $\l$. We identify a
partition with its \emph{Young diagram}, i.e., the array of left-justified squares
(boxes) with $\l_1$ boxes in the first row, $\l_2$ boxes in the second row, and so
on. The rows are arranged in matrix form from top to bottom. By the box in position
$(i,j)$ we mean the box  in the $i$-th row and $j$-th column of $\l$. The
\emph{length} of $\l$, denoted $\ell(\l)$, is the number of rows in the Young
diagram. 

Given two partitions $\l$ and $\m$, we write $\m\subseteq \l$ if and only
if $\ell(\m) \leq \ell(\l)$ and $\m_i\leq \l_i$ for $1\leq i\leq \ell(\m)$. If $\m
\subseteq \l$, we denote by $\l/\m$ the skew shape obtained by removing the boxes
corresponding to $\m$ from $\l$. The length and parts of a skew diagram are defined
in the same way as for Young diagrams.

Let $D=\l/\m$ be a skew shape and let $a=(a_1,a_2,\cdots,a_k)$ be a
sequence of positive integers such that $\sum a_i =|D|=|\l|-|\m|$. A
\emph{decomposition} of $D$ of type $a$, denoted $D_1+\cdots +D_k =D$, is given by a
sequence of shapes
$$\m=\l^{(0)} \subseteq \l^{(1)}\ldots \subseteq \l^{(k)} =\l,$$
where  $D_i=\l^{(i)}/\l^{(i-1)}$ and $|D_i|=a_i$.

A
\emph{semi-standard Young tableau} (SSYT) \emph{of shape} $\l/\m$ is a filling of
the boxes of the skew shape $\l/\m$ with positive integers so that the numbers
weakly increase in each row from left to right and strictly increase in each column
from top to bottom. The \emph{type} of a SSYT $T$ is the sequence of non-negative
integers $(t_1,t_2,\ldots)$, where $t_i$ is the number of $i$'s in $T$.

\vskip10pt
\vbox{\centering
{\young(:::225789,:4579,28,8)}
\vskip 6pt
{\scriptsize A SSYT of shape $\l=(9,5,2,1)/(3,1)$ and type
$(0,3,0,1,2,0,2,3,2)$}
 \vskip 6pt
 \centerline{\sc Figure 1}
} 
\vskip10pt

Given $T$, a SSYT  of shape $\l/\m$ and type $(t_1,t_2,\ldots)$, we
define its \emph{weight}, denoted $w(T)$, to be the monomial obtained by replacing
each $i$ in $T$ by $x_i$ and taking the product over all boxes, i.e.,
$w(T)=x_1^{t_1}x_2^{t_2}\cdots$.
The skew Schur function $s_{\l/\m}$ is defined combinatorially
by the formal power series
$$s_{\l/\m}= \sum_T w(T),$$
where the sum runs over all SSYTs of shape $\l/\m$. To obtain the usual Schur function one sets
$\m =\emptyset$.

For any positive integer $n$, the Schur function indexed by the partition $(n)$ is called
the $n$\emph{-th homogeneous symmetric function} and will be denoted by $h_n$. That is,
$h_n:=s_{(n)}.$

For any two Young diagrams $\l$ and $\m$, we let $\l \times \m$ denote
the diagram obtained by joining the corners of the leftmost, lowest box in $\l$, i.e
the box in position $(\ell(\l),1)$, with the rightmost, highest box of $\m$, i.e.,
the box in position $(1,\m_1)$. \

\vskip10pt
\vbox{\centering {
\young(:::\hfil\hfil,:::\hfil,\hfil\hfil\hfil,\hfil\hfil)} 
\vskip 6pt
{\scriptsize $\l\times \m$
where $\l=(2,1)$ and $\m=(3,2)$.} 
\vskip 6pt 
\centerline{\sc Figure 2}
}
\vskip10pt

It
follows directly from the combinatorial definition of Schur functions  that $s_{\l
\times \m} = s_{\l}s_{\m}$. One defines similarly $A \times B$, when $A$ and $B$ are
skew shapes. 

The \emph{Littlewood-Richardson coefficients} are defined
via the Hall inner product on symmetric functions  (see \cite[ pg. 306]{s}) as
follows:
$$c_{\m \n}^{\l} :=\langle s_{\l},s_{\m}s_{\n}\rangle = \langle s_{\l/\m} ,s_\n \rangle. $$
That is, $c_{\m \n}^\l$ is the coefficient of $s_\l$ in the product $s_\m s_\n$.
The Littlewood-Richardson rule gives a combinatorial interpretation for the
coefficients $c_{\m \n}^{\l}$. Before we state the rule we recall some terminology.
A \emph{lattice permutation} is a sequence $a_1a_2\cdots a_n$ such that in any
initial factor $a_1a_2\cdots a_j$, the number of $i$'s is at least as great as the
number of $(i+1)$'s for all $i$.

The \emph{reverse reading word} of a tableau is the sequence of entries of $T$ obtained
by reading the entries from right to left and top to bottom, starting with the first
row.

The \emph{Littlewood-Richardson rule} states that the coefficient $c_{\m
\n}^{\l}$ is equal to the number of SSYTs of shape $\l/\m$ and type $\n$ whose
reverse reading word is a lattice permutation. 

\begin{example} The
coefficient of $s_{(5,4,3)}$ in $s_{(4,3,2)}s_{(2,1)}$ is 2 since there are two
Littlewood-Richardson tableaux of shape $(5,4,3)/(2,1)$ and type $(4,3,2)$:

\vskip10pt
\vbox{
\centerline{ \young(::111,:222,133)
\hspace{.3in}
\young (::111,:122,233)}
}
\end{example}

The Kronecker product of Schur functions is defined via the Frobenius
characteristic map, $ch$, from the center of the group algebra of $S_n$ to the ring
of symmetric functions. For a definition of the Frobenius map, see \cite[pg.
351]{s}. The map $ch$ is a ring homomorphism and an isometry. It is known that for
any irreducible character $\chi^\l$ of the symmetric group
$$ch(\chi^\l) = s_\l.$$
Let $\chi^\l$ and $\chi^\m$ be two irreducible characters of $S_n$. The
\emph{Kronecker product} $\chi^\l\chi^\m$ is defined for every $\sigma \in S_n$ by
$\chi^\l\chi^\m(\sigma) = \chi^\l(\sigma)\chi^\m(\sigma)$.  Then
$$\chi^\l \chi^\m =\sum_{\n\vdash n} g_{\l,\m,\n} \chi^\n.$$
Hence, using the Frobenius characteristic map, one defines the Kronecker product of
Schur functions by
$$s_\l \ast s_\m := \sum_{\n\vdash n} g_{\l, \m, \n} s_\n.$$

 Littlewood \cite {l} proved the following
identity:
$$s_{\l}s_{\m} \ast s_{\eta} =
\sum_{\gamma \vdash |\l|}\sum_{\delta\vdash |\m|}c_{\gamma\delta}^{\eta} (s_{\l}\ast
s_{\gamma})(s_{\m}\ast s_{\delta}),$$ where $c_{\gamma \delta}^{\eta}$ is the
Littlewood-Richardson coefficient. Garsia and Remmel \cite{gr} generalized this result for any skew
shapes $A$, $B$ and $D$:
$$(s_A s_B)\ast s_D=
\sum_{\stackrel{D_1+D_2=D}{|D_1|=|A|, |D_2|=|B|}} (s_A\ast s_{D_1})(s_B\ast
s_{D_2}),$$ where the sum runs over all
decompositions of the skew shape $D$. By induction one obtains:
$$(h_{n_1} h_{n_2}\cdots h_{n_k})\ast s_D=
\sum_{\stackrel{D_1+\cdots +D_k=D}{|D_i|=n_i}} s_{D_1}\cdots s_{D_k},$$ where the
sum runs over all decompositions of $D$ of type $(n_1,n_2,\ldots, n_k)$. 


\section{A Schur Positivity Theorem}

In this section we consider the Schur positivity of the symmetric function
$s_{\l/\a}s_{\a}-s_{\l/\a^-}s_{\a^-}$, where $\a=(\a_1, \ldots,
\a_{\ell(\a)})\subseteq \l$ with $\a_1>\a_2$ and
$\a^-=(\a_1-1,\a_2,\ldots,\a_{\ell(\a)})$. More explicitly, we show that this
symmetric function is Schur positive if and only if $\l_1\geq 2\a_1-1$. In order to
prove this result, we need a variation of the Remmel--Whitney algorithm for expanding
the skew Schur function $s_{\l/\m}$. Recall that the \emph{reverse lexicographic filling}
of $\m$, $rl(\m)$, is a filling of the Young diagram $\m$ with the numbers $1,2, \ldots, |\m|$ 
so that the numbers are entered in order from right to left and top to bottom.

\medskip
\noindent {\bf Skew Algorithm:}
The algorithm for computing
 $$s_{\l/\m} = \sum_{|\nu|=|\l|-|\m|} c_{\m\,\n}^\l s_\n$$ is as follows:
\begin{enumerate}
\item[(1)] Form the reverse lexicographic filling of $\mu$, $rl(\m)$.
\item[(2)] Starting with the Young diagram $\lambda$,  label $|\m|$ of its outermost boxes with
the numbers $|\mu|, |\mu|-1, \ldots, 2, 1$, starting with $|\m|$, so that the
following conditions are satisfied:
\begin{enumerate}
\item[(a)] After labelling each box, the unlabelled boxes form  a Young diagram.
\item[(b)] Suppose that in $rl(\m)$ the box in
position $(i,j)$ has label $x$, where $x\leq|\m|$. If $j>1$, let $x^-$ be the label in
position $(i,j-1)$ in $rl(\m)$. If $i<\ell(\m)$, let $x^+$ be the label in position
$(i+1,j)$ in $rl(\m)$. Then in $\l$, $x$ will be placed to the left and weakly below
(to the SW) of $x^-$ and above and weakly to the right (to the NE) of $x^+$.
\end{enumerate}
\item[(3)] From each of the diagrams obtained (with $|\m|$ labelled boxes), remove
all labelled boxes. The resulting unlabelled diagrams correspond to  the summands in
the decomposition of $s_{\l/\m}$.
\end{enumerate}

\begin{remark}  Suppose $(i,j)$ is the position of $x$ in
$rl(\m)$ and  $(l,m)$ is the new position of $x$ in $\l$. The conditions (a) and (b)
impose constraints on $l$ and $m$. It can be easily verified that $l \ge i$ and $m
\ge \m_i-j+1$, where $\m_i$ is the number of boxes in the $i$-th row
of $\m$. 
\end{remark}

\begin{example} The decomposition of $s_{\l/\m}$, where
$\l=(4,4,2,2)$, $\m=(3,3)$: 
\vskip10pt

\centering{
{\Yvcentermath1 $\l=\ ${
$\yng(4,4,2,2)$}},\quad{\Yvcentermath1 $rl(\m)=$ {
      $\young(321,654)$}}.}
\vskip10pt

\noindent First we  establish the constraints on the position of each
 label in $\l$ using the Remark.

\vskip10pt
 \noindent \begin{tabular}{|c|c|l|l|}\hline label
 & position  $(i,j)$ in $\m$ & position  $(l,m)$ in $\l$ & position relative to $x^-$ and
 $x^+$\\ \hline

 $6$ & $(2,1)$ & $l \ge 2$ and $m \ge 3-1+1=3$ &\\
$5$ & $(2,2)$ & $l \ge 2$ and $m \ge 3-2+1=2$ & SW of $6$\\
  $4$ & $(2,3)$ & $l \ge 2$ and $m \ge 3-3+1=1$ & SW of $5$\\
  $3$ & $(1,1)$ & $l \ge 1$ and $m \ge 3-1+1=3$ & NE of $6$\\
$2$ & $(1,2)$ & $l \ge 1$ and $m \ge 3-2+1=2$ & SW of $3$
 and NE of $5$\\
 $1$ & $(1,3)$ & $l \ge 1$ and $m \ge 3-3+1=1$ & SW of $2$
 and NE of $4$\\ \hline
\end{tabular}
\vskip10pt

The algorithm is carried out in Fig.~3.
Thus, $s_{\l/\m}=s_{(2,2,1,1)}+ s_{(3,2,1)} + s_{(3,3)}$. 
\end{example}

The Skew Algorithm above follows from the Remmel--Whitney  algorithm \cite{rwy} and
the fact that skewing is the adjoint operation of multiplication, i.e.,
$\langle s_{\l/\m},s_{\n}\rangle=\langle s_{\l},s_{\m}s_{\n}\rangle$. In some sense we are reversing the steps
taken in \cite{rwy} when expanding the product $s_{\m}s_{\n}$. 

In order
to state our first result we need to recall the definition of \emph{lexicographic
order} on the set of all partitions. If $\l=(\l_1, \l_2, \ldots, \l_{\ell(\l)})
\vdash n$ and $\m=(\m_1, \m_2, \ldots, \m_{\ell(\m)}) \vdash m$, we say that $\l$ is
less than $\m$   \emph{in lexicographic order}, and write $\l <_{l} \m$, if there is
a non-negative integer $k$ such that $\l_i=\m_i$ for all $i=1,2, \ldots, k$ and
$\l_{k+1} <\m_{k+1}$. The lexicographic order is a total order on the set of all
partitions. 

\begin{lemma} Consider the partitions  $\l=(\l_1, \l_2, \ldots, \l_{\ell(\l)})$ and
$\a=(\a_1, \a_2, \ldots, \a_{\ell(\a)})$ such that $\a \subseteq \l$. The smallest
partition $\n$ in lexicographic order such that $s_{\n}$ appears in the expansion of
$s_{\l/\a}$ is the partition obtained by reordering the parts of $\l/\a$ in
decreasing order, i.e., the parts of $\n$ are $\l_1-\a_1, \l_2, -\a_2, \ldots,
\l_{\ell(\l)}-\a_{\ell(\l)}$ ($\a_i=0$ if $i > \ell(\a)$) reordered such that $\n$
is a partition. Moreover, the multiplicity of $s_{\n}$ in the expansion of
$s_{\l/\a}$ is equal to $1$.

\end{lemma}

\begin{proof} Using the Skew Algorithm, we obtain the smallest partition in
lexicographic order when the labels in  $rl(\a)$  are each placed in the highest possible 
row of $\l$ (since we are removing the largest possible number of boxes from the top rows of $\l$). 
We will show that the partition obtained in this way is precisely the partition obtained by
reordering the rows of $\l/\a$. We argue inductively by the number of rows of $\a$.

Assume $\a=(\a_1)$. We form the reverse lexicographic order of $\a$.
According to the Skew Algorithm, the label $|\a|$ can be placed in the first row of
$\l$ if $\l_1>\l_2$. In general, the highest position where we can place  $|\a|$  is
$(k, \l_k)$, where $k$ is the positive integer such that $\l_k=\l_1$ and
$\l_{k+1}<\l_1$. We will place the other labels of $\a$ to the SW of this position
respecting the rules of the algorithm. We will remove the highest possible
horizontal strip (a skew shape so that no two boxes are in the same column) with
$\a_1$
 \newpage

\vspace*{-2.5cm}

\begin{center}
 \includegraphics[width=6truein,height=8truein]{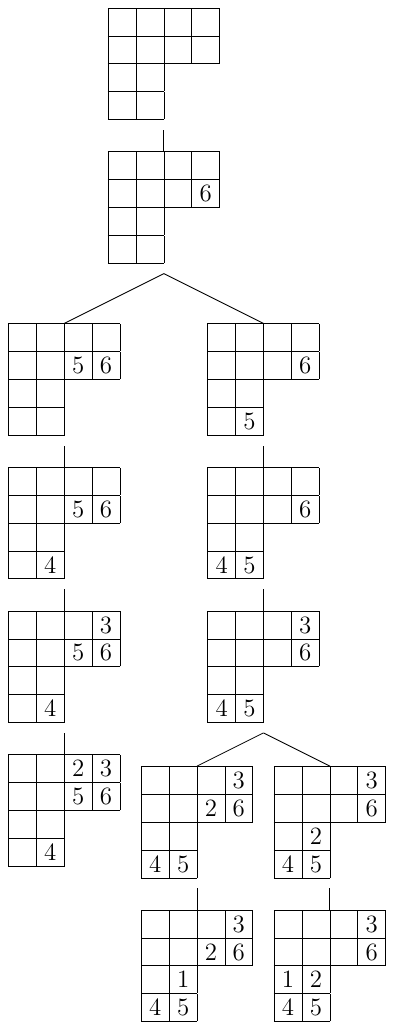}

 \end{center}\vspace*{-5cm}
 
  \centerline{\sc Figure 3}\medskip
 
\noindent boxes starting with position $(k, \l_k)$ and continuing SW. This also follows
from Pieri's rule \cite[Corollary~7.15.9]{s}. 

Now suppose $t$ is the
positive integer such that $\l_k-\l_{k+t} \geq \a_1$ and $\l_k - \l_{k+t-1}<\a_1$
(i.e., label $1$ will be placed in the $(k+t)$-th row). Then the smallest shape in
lexicographic order appearing in the Skew Algorithm is
\begin{multline*}
 (\l_1, \ldots, \l_{k-1}, \l_k-(\l_k-\l_{k+1}), \l_{k+1}-(\l_{k+1}-\l_{k+2}),
\ldots, \l_{k+t-1}-(\l_{k+t-1}-\l_{k+t}), \\ 
\kern6cm \l_{k+t}-(\a_1-(\l_k-\l_{k+t})),
\l_{k+t+1}, \ldots, \l_{\ell(\l)})\\ =(\l_1, \ldots, \l_{k-1}, \l_{k+1}, \l_{k+2},
\ldots, \l_{k+t}, \l_k-\a_1, \l_{k+t+1}, \ldots, \l_{\ell(\l)}).
\end{multline*}
 Since
$\l_k=\l_1$, this is precisely the partition obtained by reordering the rows of
$\l/(\a_1)=(\l_1-\a_1, \l_2, \ldots, \l_{\ell(\l)})$. 

Suppose the lemma
is true for all partitions with $\ell$ parts that are contained in $\l$ and let
$\a=(\a_1, \a_2,\ldots, \a_{\ell}, \a_{\ell+1})$ be a partition with $\ell+1$ parts
such that $\a \subseteq \l$. Thus $\ell(\a)=\ell+1$. 

When we place the
labels from the reverse lexicographic order of $\a$ into the boxes of $\l$ according
to the Skew Algorithm such that each label is  placed in the highest possible row,
we first place the labels of the last row of $rl(\a)$. They are placed as in the
case $\a=(\a_1)$ above, starting with placing $|\a|$ in position $(k, \l_k)$, where
$k$ is the positive integer such that $\l_k=\l_{\ell(\a)}$ and
$\l_{k+1}<\l_{\ell(\a)}$ (note that $k \geq \ell(\a)$). The remaining labels of the
last row of $\a$ are placed to the SW of this position forming a horizontal strip
with boxes in the highest possible rows of $\l$. Observe that this is equivalent to
labelling the highest possible horizontal strip of $\a_{\ell(\a)}$ boxes in
$(\l_{\ell(\a)},\l_{\ell(\a)+1},\ldots, \l_{\ell(\l)})$ such that the unlabelled
boxes yield a Young diagram, and then adding back the rows $\l_1,\ldots,
\l_{\ell(\a)-1}$ above the shape $(\l_{\ell(\a)},\l_{\ell(\a)+1},\ldots,
\l_{\ell(\l)})$ (with the labelled boxes). Notice that the unlabelled boxes form the
partition which is a rearrangement of
$(\l_1,\ldots,\l_{\ell(\a)}-\a_{\ell(\a)},\ldots, \l_{\ell(\l)})$. According to the
Skew Algorithm, we now place the labels in row $\ell(\a)-1$ of $rl(\a)$. Requiring
that the labels be placed in the highest rightmost position at each step will
automatically satisfy the requirements of the Skew Algorithm. Therefore, in order to
obtain the smallest partition in lexicographic order, we just need to label the
highest horizontal strips with $\a_{\ell(\a)}$ boxes starting at row $\ell(\a)$ and
continuing SW, then the highest horizontal strip with $\a_{\ell(\a)-1}$ boxes
starting at row $\ell(\a)-1$ and continuing SW in the remaining unlabelled boxes,
and so on until we label the highest horizontal strip with $\a_1$ boxes starting at
row $1$ such that, at each step, removing all labelled boxed yields a Young diagram.
See Fig.~4 for an example of placing the labels of the last 3 rows of $rl(\a)$ in
$\l$; the horizontal strip containing $a$'s is the one corresponding to row
$\ell(\a)$, the strip containing $b$'s corresponds to row $\ell(\a)-1$ and the strip
containing $c$'s corresponds to row $\ell(\a)-2$. 

\begin{figure}\color{black}
\centerline{
\begin{picture}(300,200)
\put(0,10){\line(0,0){190}} \put(0,200){\line(1,0){300}}
\put(300,180){\line(0,0){20}} \put(250,180){\line(1,0){50}}
\put(250,150){\line(0,0){30}} \put(200,150){\line(1,0){50}}
\put(200,120){\line(0,0){30}} \put(170,120){\line(1,0){30}}
\put(170,70){\line(0,0){50}} \put(140,70){\line(1,0){30}}
\multiput(-5,120)(10,0){18}{\line(1,0){5}}
\multiput(-5,110)(10,0){18}{\line(1,0){5}} \put(-30,110){$\ell(\a)$}
\put(140,70){\line(0,0){10}} \put(140,80){\line(1,0){30}} \put(140,70){\small{$a a
\cdots$}} \put(140,80){\line(0,0){10}} \put(140,90){\line(1,0){30}}
\put(140,80){\small{$b b \cdots$}} \put(140,90){\line(0,0){10}}
\put(140,100){\line(1,0){30}} \put(140,90){\small{$c c \cdots$}}
\put(110,40){\line(0,0){10}} \put(110,50){\line(1,0){30}} \put(110,40){\small{$a a
\cdots$}} \put(80,40){\line(0,0){10}} \put(80,50){\line(1,0){30}}
\put(80,40){\small{$c c \cdots$}} \put(110,50){\line(0,0){10}}
\put(110,60){\line(1,0){30}} \put(110,50){\small{$b b \cdots$}}
\put(110,60){\line(0,0){10}} \put(110,70){\line(1,0){30}} \put(110,60){\small{$c c
\cdots$}} \put(170,120){\line(0,0){10}} \put(170,130){\line(1,0){30}}
\put(170,120){\small{$b b \cdots$}} \put(170,130){\line(0,0){10}}
\put(170,140){\line(1,0){30}} \put(170,130){\small{$c c \cdots$}}
\put(140,40){\line(0,0){30}} \put(70,40){\line(1,0){70}} \put(70,30){\line(0,0){10}}
\put(30,30){\line(1,0){40}} \put(30,10){\line(0,0){20}} \put(0,10){\line(1,0){30}}
\put(110,0){\sc Figure 4}
\end{picture}}
\end{figure}

If we follow this procedure of placing each label as high as possible we
obtain the smallest partition in lexicographic order appearing in the expansion of
$s_{\l/\a}$. The labels in rows $2,3,\ldots, \ell(\a)$ of $rl(\a)$ can only be
placed in $\bar{\l}=(\l_2,\l_3,\ldots,\l_{\ell(\l)})$. The unlabelled boxes form the
partition $\m=(\l_1,\m_2, \ldots,\m_{\ell(\m)})$ where $(\m_2,\ldots,
\m_{\ell(\m)})$ is the smallest partition in lexicographic order occurring in
$s_{\bar{\l}/\bar{\a}}$. By induction hypothesis, $(\m_2,\ldots, \m_{\ell(\m)})$ is
obtained by rearranging $\l_2-\a_2, \l_3-\a_3, \ldots,
\l_{\ell(\a)}-\a_{\ell(\a)},\l_{\ell(\a)+1}, \ldots \l_{\ell(\l)}$. If we continue
the labelling in $\l$ with the labels in row $\a_1$ of $rl(\a)$, the labels are
placed in $\m=(\l_1,\m_2,\ldots, \m_{\ell(\m)})$. By the discussion in the previous
paragraph, we will be labelling the highest possible horizontal strip in $\m$.
Hence, removing this strip from $\m$ yields the partition $\n$ which is a
rearrangement of $\l_1-\a_1,\m_2, \ldots, \m_{\ell(\m)}$. By induction we have that
$\nu$ is the smallest partition and it is a rearrangement of $\l_1-\a_1,\l_2-\a_2,
\ldots \l_{\ell(\l)}-\a_{\ell(\l)}$, where $\a_i=0$ for $i>\ell(\a)$. 

Since there is only one way of placing the labels of the reverse lexicographic order
of $\a$ in the highest possible rows of $\l$,  the multiplicity of $s_{\n}$ in
$s_{\l/\a}$, where $\n$ is the smallest partition in lexicographic order, equals
$1$. \end{proof}

\begin{cor} Let $\a=(\a_1, \a_2, \ldots, \a_{\ell(\a)})$ and $\b=(\b_1, \b_2,
\ldots, \b_{\ell(\b)})$ be arbitrary partitions. The smallest partition $\n$ in
lexicographic order appearing in the expansion of $s_{\a}s_{\b}$ is the partition
obtained by concatenating the parts of $\a$ and $\b$ and reordering them to form a
partition.
\end{cor}

\begin{proof} The proof is a direct consequence of the Littlewood-Richardson
rule. It also follows directly from the Remmel--Whitney algorithm for multiplying two
Schur functions \cite{rwy} and Lemma~2.1. \end{proof}

\begin{definition} A symmetric function is said to be \emph{Schur positive}
(or \emph{s-positive}) if, when expanded as a linear combination of Schur functions,
all the coefficients are positive.
\end{definition}

Consider the product $s_{\l/\a}s_{\a}$. The combinatorial definition of Schur
functions implies that $s_{\l/\a}s_{\a}$ is the skew Schur function corresponding to
the skew shape $\m/\eta$, where $$\m=(\l_1+\a_1, \l_1+\a_2,\ldots
\l_1+\a_{\ell(\a)}, \l_1,\l_2, \ldots, \l_{\ell(\a)}, \l_{\ell(\a)+1}, \ldots ,
\l_{\ell(\l)})$$ and $$\eta=(\underbrace{\l_1, \ldots , \l_1}_{\ell(\a) \
\mbox{times}}, \a_1, \a_2, \ldots, \a_{\ell(\a)}).$$ This is the skew shape
$\a\times \l/\a$.

\begin{example} Let
$\l=(6,4,2,2)$ and $\a=(3,1)$. Then $s_{\l/\a}s_{\a}$ is the skew Schur function
corresponding to the skew shape $\a \times \l/\a$ below.

\vskip10pt
\vbox{
{\centerline {
\young(::::::\hfil\hfil\hfil,::::::\hfil,:::\hfil\hfil\hfil,:\hfil\hfil\hfil,\hfil\hfil,\hfil\hfil)}}
}
\end{example}

\begin{definition} Let $\a=(\a_1,\a_2,\, \ldots, \a_{\ell(\a)})$ be any
sequence of non-negative integers. A sequence $a_1a_2 \cdots a_n$ is an
\emph{$\a$-lattice permutation} if in any initial factor  $a_1a_2\cdots a_j$, $1\leq
j\leq n$, we have for any positive integer $i$:
$$\mbox{the number of } i's  +\a_i \geq \mbox{ the number of }
(i+1)'s + \a_{i+1}.$$ Here $\a_i=0$ if $i >\ell(\a)$. 
\end{definition}

Then, if
$\n=(\n_1, \n_2, \ldots, \n_{\ell(\n)})\vdash n$, the multiplicity of $s_{\n}$ in
$s_{\l/\a}s_{\a}$ is given by the number of SSYT of shape $\l/\a$ and type
$\n/\a:=(\n_1-\a_1, \n_2-\a_2, \ldots, \n_{\ell(\a)}-\a_{\ell(\a)}, \n_{\ell(\a)+1},
\ldots, \n_{\ell(\n)})$ whose reverse reading word is an $\a$-lattice permutation.
(If $\a \not \subseteq \n$, the multiplicity of $s_{\n}$ in $s_{\l/\a}s_{\a}$ is
$0$.) 

\begin{theorem} Let $\l\vdash n $, $\a=(\a_1,\a_2,\ldots \a_{\ell(\a)})\vdash p$
with $\a_1>\a_2$ and $\a^-=(\a_1-1,\a_2,\ldots \a_{\ell(\a)})\vdash p-1$. Assume
that $\a \subseteq \l$. The symmetric function $s_{\l/\a}s_{\a}-
s_{\l/\a^-}s_{\a^-}$ is Schur positive if and only if $\l_1 \geq 2 \a_1-1$.
\end{theorem}

\begin{proof} Assume that $\l_1 < 2\a_1-1$. Let $\n$ denote the smallest partition in
lexicographic order such that $s_{\n}$ appears in the expansion of $s_{\l/\a}s_{\a}$
and let $\n^*$ denote the smallest partition in lexicographic order such that
$s_{\n^*}$ appears in the expansion of $s_{\l/\a^-}s_{\a^-}$. Since
$s_{\l/\a}s_{\a}$ is the skew Schur function corresponding to the skew shape
$\a\times \l/\a$, it follows from  Lemma~2.1 that the parts of $\n$ are precisely
the parts of $\a\times \l/\a$ reordered to form a partition. Since $\a_i=0$ for $i
>\ell(\a)$, the parts of $\n$ are $ \a_1, \a_2, \ldots,\a_{\ell(\a)}, \l_1-\a_1,
\l_2-\a_2,\ldots \l_{\ell(\a)}-\a_{\ell(\a)},\l_{\ell(\a)+1}, \ldots,
\l_{\ell(\l)}$. Similarly, the parts of $\n^*$ are precisely the parts of
$\a^-\times \l/\a^-$, that is $ \a_1-1, \a_2, \ldots, \a_{\ell(\a)}, \l_1-\a_1+1,
\l_2-\a_2, \ldots,\l_{\ell(\a)}-\a_{\ell(\a)}, \l_{\ell(\a)+1}, \ldots,
\l_{\ell(\l)}$. The partitions $\n$ and $\n^*$ differ in only two parts: $\l_1-\a_1$
vs. $\l_1-\a_1+1$ and $\a_1$ vs. $\a_1-1$. Since $\l_1<2\a_1-1$, we have
$\l_1-\a_1<\a_1-1$ and $\l_1-\a_1+1<\a_1$. Thus, $\l_1-\a_1<\a_1$ and
$\l_1-\a_1+1\leq\a_1-1$. Let $k$ be the positive integer such that $\n_k=\a_1$ and
$\n_{k+1}<\a_1$. Then $\n^*_k=\a_1-1$. We have $\n_i=\n^*_i$ for all $i=1,2, \ldots,
k-1$ and $\n^*_k<\n_k$. Thus $\n^*$ is lexicographically smaller than $\n$ and
$s_{\n^*}$ appears in the expansion of $s_{\l/\a}s_{\a}- s_{\l/\a^-}s_{\a^-}$ with
coefficient $-1$. 

Assume now that $\l_1\geq 2\a_1-1$. We will show that,
for any $\n \vdash n$,  the coefficient of $s_{\n}$ in $s_{\l/\a^-}s_{\a^-}$ is
smaller or equal than its coefficient in $s_{\l/\a}s_{\a}$. Let $\n \vdash n$ be
such that $s_{\n}$ appears in the expansion of $s_{\l/\a^-}s_{\a^-}$. The
coefficient of $s_{\n}$ in the expansion of $s_{\l/\a^-}s_{\a^-}$ is equal to the
number of SSYTs of shape $\l/\a^-$ and type $\n/\a^-$ whose reverse reading word is
an $\a^-$-lattice permutation. Let $T$ be a SSYT of shape $\l/\a^-$ and type
$\n/\a^-$ whose reverse reading word is an $\a^-$-lattice permutation. Then, in the
first row of $\l/\a^-$,
$$\begin{array}{rll} \mbox{ the
number of } (i+1)'s + \a_{i+1}\kern-6pt &\leq \a_i,& \mbox{ if } i \geq 2 \mbox { and } \\
\mbox{ the number of } 2's + \a_2 \kern-6pt&\leq \a_1-1. \end{array}$$
Recall that $\a_j=0$ if $j >\ell(\a)$. 

Thus, the  number of boxes in the
first row of $\l/\a^-$ with labels different from $1$ is at most
$$\sum_{i\geq 2} (\a_i-\a_{i+1}) + \a_1-1-\a_2= \a_1-1.$$
Since $\l_1 \geq 2 \a_1-1$, the number of boxes in the first row of $\l/\a^-$ satisfies
$$\l_1-\a_1+1\geq \a_1>\a_1-1.$$
Thus, in $T$,  the leftmost box  in the first row of $\l/\a^-$,  $(1,\a_1)$,  must be filled with $1$.
 
Let $T^*$ be the tableau
obtained from $T$ by removing the leftmost box in the first row of $T$, i.e., the box
labelled $1$ in position $(1, \a_1)$. Then $T^*$ is a SSYT of shape $\l/\a$ and type
$\n/\a$ whose reverse reading word is an $\a$-lattice permutation. 

For
every SSYT of shape $\l/\a^-$ and type $\n/\a^-$ whose reverse reading word is an
$\a^-$-lattice permutation we obtained a SSYT of shape $\l/\a$ and type $\n/\a$
whose reverse reading word is an $\a$-lattice permutation. Thus the multiplicity of
$s_{\n}$ in $s_{\l/\a}s_{\a}$ is greater or equal to the multiplicity of $s_{\n}$ in
$s_{\l/\a^-}s_{\a^-}$ and $s_{\l/\a}s_{\a}- s_{\l/\a^-}s_{\a^-}$ is
Schur-positive.\end{proof}


\section{Combinatorial rule for the Kronecker coefficient}

In this section we use Theorem~2.3 to give a combinatorial rule for computing the
coefficient $g_{(n-p,p),\l,\n}$  in the Kronecker product
$s_{(n-p,p)}\ast s_{\l}$, whenever $\l_1\geq 2p-1$ or $\ell(\l) \geq 2p-1$.  Our 
combinatorial rule works for all partitions $\l$ that do not fit in the $(2p-2)
\times (2p-2)$ square. If $n> (2p-2)^2$, then the rule applies to all partitions
$\l$ of $n$.

\begin{definition} A SSYT $T$ of shape $\l/\a$
and type $\n/\a$ whose reverse reading word is an $\a$-lattice permutation is called
a \emph{Kronecker Tableau} of shape $\l/\a$ and type $\n/\a$ if
\begin{itemize}
\item[(I)] $\a_1=\a_2$ or
\item[(II)] $\a_1>\a_2$ and  any one of the following two conditions is
satisfied:
\begin{itemize}
\item[(i)]  The number of $1$'s in the second row of $\l/\a$ is exactly $\a_1-\a_2$.
\item[(ii)]  The number of $2$'s in the first row of $\l/\a$ is exactly $\a_1-\a_2$.
\end{itemize}
\end{itemize}

\noindent
Denote by $k_{\a\n}^{\l}$ the number of Kronecker tableaux of
shape $\l/\a$ and type $\n/\a$. 
\end{definition}

\begin{note} There are $\a_1-\a_2$ 1's 
in the second row of $\l/\a$ if and only if there is a $1$ in box $(2,\a_1)$ of $\l/\a$.
\end{note}
\begin{lemma}
Let $\l\vdash n$ and $\a\vdash p$ be partitions such that $\a_1>\a_2$, $\a\subseteq \l$
and $\l_1\geq 2\a_1-1$ and let $\a^-=(\a_1-1,\a_2,\ldots, \a_{\ell(\a)})$. Then the
coefficient of $s_{\n}$ in $s_{\a}s_{\l/\a}-s_{\a^-}s_{\l/\a^-}$ equals the number
of Kronecker tableaux of shape $\l/\a$ and type $\n/\a$, i.e., $k^{\l}_{\a\n}$.
\end{lemma}
\begin{proof} By Theorem~2.3 we have that $s_{\a}s_{\l/\a}-s_{\a^-}s_{\l/\a^-}$ is Schur
positive if and only if $\l_1\geq 2\a_1-1$. Furthermore, we have observed that the
coefficient of $s_\n$ in $s_{\a}s_{\l/\a}$ is equal to the number of SSYTs of shape
$\l/\a$ and type $\n/\a$ whose reverse reading word is an $\a$-lattice permutation.
Denote the set of these SSYTs by $\mathcal{T}_{\l/\a}^\n$ and let $K_{\l/\a}^{\n}$
denote the set of Kronecker tableaux of shape $\l/\a$ and type $\n/\a$. Denote by
$NK_{\l/\a}^{\n}$ the set of elements in $\mathcal{T}_{\l/\a}^\n$ that are not
Kronecker tableaux. Hence, $\mathcal{T}_{\l/\a}^\n=K_{\l/\a}^{\n} \cup
NK_{\l/\a}^{\n}$ (disjoint union). 

Let $\a$ be such that $\a_1>\a_2$ and
let $\mathcal{T}_{\l/\a^-}^\n$ be the set of SSYTs of shape $\l/\a^-$ and type
$\n/\a^-$ whose reverse reading word is an $\a^-$-lattice permutation. We show that
there exists a bijection
$$\mathcal{T}_{\l/\a^-}^\n\longrightarrow NK_{\l/\a}^\n.$$
In Theorem~2.3, we showed that in every $T\in \mathcal{T}_{\l/\a^-}^\n$ the box in
position  $(1,\a_1)$ is filled with $1$. The bijection is given by removing this box
from $T$ to obtain a tableau $T^\ast \in \mathcal{T}_{\l/\a}^\n$. Furthermore, the
number of 2's in the first row in $T^\ast$ can be at most $\a_1-\a_2-1$ since $T$
yielded an $\a^-$-lattice permutation. Also, the label in position $(2,\a_1)$ of
$T^\ast$ must be greater than 1, since in $T$ the box in position $(1, \a_1)$ was
filled with $1$. Thus $T^\ast \in NK_{\l/\a}^{\n}$. On the other hand,
  to any element $T\in
NK_{\l/\a}^{\n}$ we can add a box labelled $1$ in position $(1,\a_1)$ to obtain a tableau in
$\mathcal{T}_{\l/\a^-}^{\n}$.  Thus, the coefficient of $s_{\n}$ in
$s_{\a}s_{\l/\a}-s_{\a^-}s_{\l/\a^-}$ equals $k^{\l}_{\a\n}$. \end{proof}

Let $\l, \n\vdash n$. The \emph{intersection}  of $\l$ and $\n$, denoted $\l\cap
\n$, is the Young diagram consisting of the boxes $(i,j)$ that belong to both $\l$
and $\n$. 

\begin{example} If {\Yvcentermath1 $\l=$ {\tiny
$\yng(5,3,1)$ }} and {\Yvcentermath1$\n=$ {\tiny$\yng(3,2,2,1,1)$}}, then
{\Yvcentermath1 $\l\cap \n=$ {\tiny$\yng(3,2,1)$.}}
\end{example}
\begin{theorem} Let $n$ and $p$ be positive integers such that $n \geq 2p$. Let $\l,
\n \vdash n$, $\l=(\l_1, \ldots, \l_{\ell(\l)})$ and let $\l'$ denote the partition
conjugate to $\l$.
\begin{itemize}\item[(a)] If $\l_1 \geq 2p-1$,   the multiplicity of $s_{\n}$ in $
s_{(n-p,p)}*s_{\l}$ equals ${\displaystyle \sum_{\stackrel{\a \vdash p}{\a \subseteq
\l\cap \nu}}k_{\a\n}^{\l}.}$ \item[(b)] If $\ell(\l) \geq 2p-1$, the multiplicity of
$s_{\n}$ in $ s_{(n-p,p)}*s_{\l}$ equals ${\displaystyle \sum_{\stackrel{\a \vdash
p}{\a \subseteq \l'\cap \nu'}}k_{\a\n'}^{\l'}.}$\end{itemize}
\end{theorem}
\begin{proof}
(a) From the Garsia--Remmel formula it easily follows that \[s_{(n-p,p)}*s_{\l}=
\sum_{\stackrel{\a \vdash p}{\a \subseteq \l}} s_{\a}s_{\l/\a} - \sum_{\stackrel{\b
\vdash p-1}{\b \subseteq \l}} s_{\b}s_{\l/\b}.
\]
There is a $1-1$ correspondence between the partitions $\a \vdash p$ with
$\a_1>\a_2$ and the partitions $\a^- \vdash p-1$ given by \[\a=(\a_1,\a_2, \ldots
\a_{\ell(\a)}) \longleftrightarrow \a^-=(\a_1-1,\a_2, \ldots \a_{\ell(\a)}).\] The
Young diagram of $\a^{-}$ is obtained from the Young diagram of $\a$ by removing the
last box of the first row. Hence, we have
\[s_{(n-p,p)}*s_{\l}= \sum_{\stackrel{\a \vdash p, \a_1>\a_2}{\a \subseteq \l}}
(s_{\a}s_{\l/\a} -  s_{\a^-}s_{\l/\a^-}) + \sum_{\stackrel{\a \vdash p,
\a_1=\a_2}{\a \subseteq \l}} s_{\a}s_{\l/\a}.
\]
 If $\a_1=\a_2$, by definition,  every SSYT whose reverse reading word is an $\a$-lattice
permutation is a Kronecker tableau. If $\a_1>\a_2$,  by Lemma~3.1  the coefficient
of $s_{\n}$ in $s_{\a}s_{\l/\a}-s_{\a^-}s_{\l/\a^-}$ equals $k^{\l}_{\a\n}$. Note
that since $\a_1 \leq p$, we have $\l_1 \geq 2 \a_1-1$ for all partitions $\a \vdash
p$.

\medskip
\noindent (b) If $\ell(\l)\geq 2p-1$, then the first part of  $\l'$ is greater than
or equal to $2p-1$.  The identity $s_{\l'}=s_{(1^n)}\ast s_{\l}$ and the
commutativity and associativity of the Kronecker product  give $s_{(n-p,p)}\ast
s_{\l}=(s_{(n-p,p)} \ast s_{\l'})*s_{(1^n)}=\left( \sum
g_{(n-p,p),\l',\n'}s_{\n'}\right)*s_{(1^n)}=\sum
g_{(n-p,p),\l',\n'}(s_{\n'}*s_{(1^n)})= \sum g_{(n-p,p),\l',\n'}s_{\n}.$ Part (a) of
the theorem applied to $\l'$ implies the result of part (b).
\end{proof} 
\noindent {\bf Example:} If $n=15$, $p=3$, $\l=(6,4,4,1)$ and
$\n=(5,4,3,3)$, the Kronecker tableaux are:

\vskip10pt
\centerline{
{ \young(:::222,1133,2344,4)} \hskip .1in {
\young(::1112,:233,2344,4)} \hskip .1in { \young(::1123,:123,2344,4)}
\hskip .1in { \young(::1123,:223,1344,4)}} 
\vskip10pt

\noindent Hence, $g_{(12,3),(6,4,4,1),(5,4,3,3)}=4$. 

\begin{note} If $s_{\nu}$ appears in $ s_{(n-p,p)}*s_{\l}$, the length of
$\nu$ satisfies $$\ell(\nu) \leq \ell(\l) + \min\{p, \ell(\l)\}.$$ In particular we obtain the 
known fact that, if $s_{\nu}$ appears in 
$s_{(\l_1,\l_2)}*s_{(\m_1,\m_2)}$, then $\nu$ has at most four parts.
\end{note}

\begin{cor}
Let $\lambda, \nu \vdash n$, then for a positive integer $p$
$$g_{(n-p,p), \l ,\n} \leq \sum_{\stackrel{\a \vdash p}{\a
\subseteq \l\cap\n}}k_{\a\n}^{\l}.$$
\end{cor}

This corollary is a direct consequence of the proof of Theorem~3.2. 

\begin{example} Consider the Kronecker product $s_{(4,3)}\ast s_{(4,3)}$.
Using MAPLE one checks that the coefficient of each $s_{(4,2,1)}$ and $s_{(3,2,2)}$
in $s_{(4,3)}\ast s_{(4,3)}$ is 1. However, there are two Kronecker tableaux of
shape $(4,3)/(2,1)$ and type $(4,2,1)/(2,1)$ and two Kronecker tableaux of shape
$(4,3)/(2,1)$ and type $(3,2,2)/(2,1)$. These are the only two Schur functions
appearing in the decomposition of $s_{(4,3)}\ast s_{(4,3)}$ with coefficient
strictly less than the number of Kronecker tableaux.
\end{example}

\section{Applications of the Combinatorial Rule}
In this section we illustrate the usefulness of Theorem~3.2. We apply the theorem to
address the question of characterizing multiplicity free Kronecker products when one
partition has two parts. We also consider special cases for $\l$ and $\n$ and obtain
simple formulas for $g_{(n-p,p),\l,\n}$ just by counting Kronecker tableaux. We believe
that, when compared with known results \cite{r1}, \cite{r2}, \cite{rwh}, \cite{ro}, the formulas
of Theorem~3.2 provides elegant combinatorial solutions to the Kronecker product coefficients
for the cases in which it works.

 \subsection{Multiplicity-Free Kronecker Products}

In this subsection we use our main theorem to determine the partitions $\lambda$ for
which all $s_\n$ in the decomposition of $s_{(n-p,p)}\ast s_{\lambda}$ have coefficients either 0 or 1, i.e., 
$s_{(n-p,p)}\ast s_{\lambda}$ is multiplicity-free.  Recall that a
partition is called \emph{rectangular} or a \emph{rectangle partition} if all its
rows have the same size. That is, $\l=(a^k) = (a, a, \ldots, a)$ is a rectangle.
Define $C(\l) := |\{i \, |\, \l_i>\l_{i+1}, 1\leq i \leq \ell(\l)-1\}|$.

\begin{prop}
\label{P:p=1}
Let $n$ be a positive integer and $\l\vdash n$. Then $s_{(n-1,1)}\ast s_\l$ is
multiplicity free if and only if
$\l$ is a rectangular partition, $\l=(a^k)$, or $\l=(a^{k_1}, b^{k_2})$.
\end{prop}
\begin{proof} It is well-know  that $s_{(n-1,1)}\ast s_\l=s_{(1)} s_{\l/(1)}-s_\l$, see for
example \cite[Exercise 7.81]{s}.  From this and  Pieri's rule we can easily see that
$$s_{(n-1,1)}\ast s_{\l} = C(\l) s_{\l} +\sum s_\m,$$

\noindent where $C(\l)$ is as defined above and  the sum
is over all partitions different from $\l$ that can be obtained by removing one box
and then adding a box to $\l$.   Hence, the only partition that occurs with
multiplicity greater than 1 is $\l$ itself when $C(\l)\geq 1$.
\end{proof}
\begin{cor}
\label{cor:cor2} Let $n$ be a positive integer and $\l\vdash n$. Then
$s_{(2,1^{n-2})}\ast s_\l$ is multiplicity free if and only if $\l$ is a rectangular
partition, $\l=(a^k)$, or $\l=(a^{k_1}, b^{k_2})$.
\end{cor}
\begin{proof} This follows from the identity $s_{(1^n)}\ast s_{(n-1,1)} =s_{(2,1^{n-2})}$
and Proposition~\ref{P:p=1}.\end{proof}

\begin{lemma}
\label{Lemma:free} Let $n, t$ be positive integers and $\lambda\vdash n$. If
$\ell(\l)\geq t$, then $s_{(1^t)}s_{\lambda/(1^t)}$ is multiplicity free if and only
if $\ell(\l)= t$ or $\lambda=(a^k)$ for some positive integers $a, k$.
\end{lemma}

\begin{proof} We first show that $s_{(1^t)}s_{\lambda/(1^t)}$ is not multiplicity free if
$\ell(\l)>t$ and $\l$ is not a rectangle partition.   We will show that the
multiplicity of $s_\l$ in  $s_{(1^t)}s_{\lambda/(1^t)}$ is greater than 1.  Recall
that the multiplicity of $s_\l$ in  $s_{(1^t)}s_{\lambda/(1^t)}$ is
$k_{(1^t)\l}^{\l}$, the number of Kronecker tableaux of shape $\l/{(1^t)}$ and type
$\l/(1^t)$.  Below we obtain two different Kronecker tableaux of shape $\l/{(1^t)}$
and type $\l/(1^t)$:
\begin{enumerate}
\item[(1)] Fill row $i$ of $\l/(1^t)$ with $i$'s.
 \item[(2)]  Let $s$ be
the number of rows of size $\l_1$ in $\l$.  Note that $s< \ell(\l)$ since $\l$ is
not a rectangular partition.

\medskip
\noindent {\sc Case I.}  If $s\leq t$. In box $(s,\l_1)$ place $t+1$ and in box
$(t+1,1)$ place $s$. Then fill the remaining boxes with $i$ if the boxes are in the
$i$-th row. 

\medskip
\noindent {\sc Case II.} If $s> t$.  In boxes $(t,\l_1)$,
$(t+1,\l_1), \ldots, (s,\l_1)$ place $t+1, t+2, \ldots, s+1$ respectively. And in
boxes $(t+1, 1), (t+2, 1), \ldots, (s+1, 1)$ place $t, t+1, t+2, \ldots, s$,
respectively. Finally fill the remaining unlabelled boxes with $i$ if the box is in
the $i$-th row.
\end{enumerate}

We now show that if $\ell(\l) \leq t$, then $s_{(1^t)}s_{\lambda/(1^t)}$
is multiplicity free.  If  $\ell(\l)<t$ then $s_{(1^t)}s_{\lambda/(1^t)}=0$, hence
multiplicity free. If $\ell(\l)=t$, then $\lambda/(1^t)=(\l_1-1,\ldots, \l_t -1)$ is
a partition. Hence  $s_{(1^t)}s_{\lambda/(1^t)}$ is multiplicity free by Pieri's
rule. 

Finally, if $\l = (a^k)$, $k \geq t$, then $s_{\lambda/(1^t)} =
s_{(a^{k-t}, (a-1)^t)}$ by \cite[Theorem~2.1]{vw}. This can also be seen by applying
the Skew Algorithm of Section~2.  In this case $s_{(1^t)}s_{\lambda/(1^t)}$ is
multiplicity free by Pieri's rule.
\end{proof}

Recall that $k_{\a \n}^{\l}$ denotes the number of Kronecker tableaux of shape
$\l/\a$ and type $\n/\a$.

\begin{lemma}
\label{Lemma:2row}
Let $n\geq 6$ and $\l=(\l_1,\l_2)\vdash n$ such that $\l_1>\l_2>1$.
\begin{enumerate}
\item[(i)] If $\l_1> \l_2+1$, then $k_{(2),(\l_1-1,\l_2,1)}^{\l}\geq 1$ and
$k_{(1,1),(\l_1-1,\l_2,1)}^{\l}\geq 1$. \item[(ii)] If $\l_1=\l_2+1$, then
$k_{(2),(\l_1,\l_2-1,1)}^{\l}\geq 1$ and $k_{(1,1),(\l_1,\l_2-1,1)}^{\l}\geq 1$.
\end{enumerate}
\end{lemma}
\begin{proof} 
(i)  Suppose $\l_1> \l_2+1$, then $\l_1-\l_2\geq 2$. The
 following is a Kronecker tableaux of shape $\l/(2)$ and type $(\l_1-1,\l_2,
1)/(2)$: 

\vskip10pt
\begin{picture}(400,25)(0,0)
\put(180,0){\line(1,0){80}} \put(200,20){\line(1,0){120}}
\put(180,10){\line(1,0){140}} \put(320,10){\line(0,0){10}}
\put(180,0){\line(0,0){10}}\put(190,0){\line(0,0){10}}\put(183,2){\scriptsize{1}}\put(193,2){\scriptsize{2}}
\put(200,0){\line(0,0){20}}\put(210,0){\line(0,0){20}}\put(203,2){\scriptsize{2}}\put(218,2){$\cdots$}
\put(240,0){\line(0,0){20}}\put(250,0){\line(0,0){20}}\put(243,2){\scriptsize{2}}\put(253,2){\scriptsize{3}}
\put(260,0){\line(0,0){20}}\put(290,10){\line(0,0){10}}\put(203,12){\scriptsize{1}}\put(218,12){$\cdots$}
\put(300,10){\line(0,0){10}}\put(310,10){\line(0,0){10}}\put(243,12){\scriptsize{1}}
\put(253,12){\scriptsize{1}}\put(269,12){$\cdots$}\put(293,12){\scriptsize{1}}
\put(303,12){\scriptsize{2}}\put(313,12){\scriptsize{2}}
\put(180,10){\red{\framebox(20,10)[t]{}}}\multiput(180,10)(.5,0){40}{\red{\line(0,1){10}}}
\end{picture}
\vskip10pt

Hence, $k_{(2),(\l_1-1,\l_2,1)}^{\l}\geq 1$. The following is a
Kronecker tableaux of shape $\l/(1,1)$ and type
$(\l_1-1,\l_2,1)/(1,1)$: 

\vskip10pt
\begin{picture}(400,25)(0,0)
\put(180,0){\line(1,0){60}} \put(180,20){\line(1,0){120}}
\put(180,10){\line(1,0){120}}
\put(180,0){\line(0,0){10}}\put(190,0){\line(0,0){20}}\put(183,2){\scriptsize{2}}\put(183,12){\scriptsize{1}}
\put(200,0){\line(0,0){20}}\put(193,2){\scriptsize{2}}\put(193,12){\scriptsize{1}}\put(208,2){$\cdots$}
\put(240,0){\line(0,0){20}}\put(208,12){$\cdots$}\put(230,0){\line(0,0){20}}
\put(232,2){\scriptsize{2}}\put(232,12){\scriptsize{1}}\put(253,12){$\cdots$}
\put(280,10){\line(0,0){10}}\put(290,10){\line(0,0){10}}\put(180,10){\line(0,0){10}}
\put(300,10){\line(0,0){10}}\put(283,12){\scriptsize{1}}\put(293,12){\scriptsize{3}}
\put(170,0){\red{\framebox(10,20)[t]{}}}\multiput(170,0)(0,.5){40}{\red{\line(1,0){10}}}
\end{picture}
\vskip10pt

Hence, $k_{(1,1),(\l_1-1,\l_2,1)}^{\l}\geq 1$. 

\medskip 
\noindent
(ii) Let $\l_1=\l_2+1$. The following is a Kronecker tableaux of shape $\l/(2)$ and
type $(\l_1,\l_2-1,1)/(2)$: 

\vskip10pt
\begin{picture}(400,25)(0,0)
\put(180,0){\line(1,0){80}} \put(200,20){\line(1,0){70}}
\put(180,10){\line(1,0){90}} \put(270,10){\line(0,0){10}}
\put(180,0){\line(0,0){10}}\put(190,0){\line(0,0){10}}\put(183,2){\scriptsize{1}}\put(193,2){\scriptsize{1}}
\put(200,0){\line(0,0){20}}\put(210,0){\line(0,0){20}}\put(203,2){\scriptsize{2}}\put(218,2){$\cdots$}
\put(240,0){\line(0,0){20}}\put(250,0){\line(0,0){20}}\put(243,2){\scriptsize{2}}\put(253,2){\scriptsize{3}}
\put(260,0){\line(0,0){20}}\put(203,12){\scriptsize{1}}\put(218,12){$\cdots$}\put(243,12){\scriptsize{1}}
\put(253,12){\scriptsize{2}}\put(263,12){\scriptsize{2}}
\put(180,10){\red{\framebox(20,10)[t]{}}}\multiput(180,10)(.5,0){40}{\red{\line(0,1){10}}}
\end{picture}
\vskip10pt

Thus, $k_{(2),(\l_1,\l_2-1,1)}^{\l}\geq 1$. The following is a Kronecker
tableaux of shape $\l/(1,1)$ and type $(\l_1,\l_2-1,1)/(1,1)$: 

\vskip10pt
\begin{picture}(400,25)(0,0)
\put(180,0){\line(1,0){70}} \put(180,20){\line(1,0){80}}
\put(180,10){\line(1,0){80}}
\put(180,0){\line(0,0){10}}\put(190,0){\line(0,0){20}}\put(183,2){\scriptsize{2}}\put(183,12){\scriptsize{1}}
\put(200,0){\line(0,0){20}}\put(193,2){\scriptsize{2}}\put(193,12){\scriptsize{1}}\put(208,2){$\cdots$}
\put(240,0){\line(0,0){20}}\put(208,12){$\cdots$}\put(230,0){\line(0,0){20}}
\put(232,2){\scriptsize{2}}\put(232,12){\scriptsize{1}}\put(180,10){\line(0,0){10}}
\put(250,0){\line(0,0){20}}\put(260,10){\line(0,0){10}}\put(243,2){\scriptsize{3}}
\put(243,12){\scriptsize{1}}\put(253,12){\scriptsize{1}}
\put(170,0){\red{\framebox(10,20)[t]{}}}\multiput(170,0)(0,.5){40}{\red{\line(1,0){10}}}
\end{picture}
\vskip10pt

Therefore, $k_{(1,1),(\l_1,\l_2-1,1)}^{\l}\geq 1$.\end{proof}

\begin{theorem}
\label{T:p=2} Let $n\geq 6$ and $\l\vdash n$.  Then $s_{(n-2,2)}\ast s_\l$ is
multiplicity free if and only if $\l=(n-1,1), (n), (1^n), (2, 1^{n-2})$ or $(m^k)$
for $m, k$ positive integers.
\end{theorem}
\begin{proof} Since $n\geq 6$, we  have  $\l_1\geq 3$ or $\ell(\l)\geq 3$. We assume
without loss of generality that $\l_1\geq 3$, since if $\l_1< 3$, then $\ell(\l)\geq
3$ and in this case  the first row of $\l'$ is greater or equal to 3 and
$s_{(1^n)}\ast (s_{(n-2,2)}\ast s_{\l'})=s_{(n-2,2)}\ast s_\l$. 

By
Theorem~3.2, $g_{(n-2,2),\l,\n} =k_{(2),\n}^\l +
k_{(1,1),\n}^{\l}$. By Lemma~\ref{Lemma:free}, 
if $\ell(\l)>2$ or if $\l$ is not a rectangular partition then
there exists a $\n\vdash n$ such that $k_{(1,1),\n}^\l\geq 2$. Therefore,
$s_{(n-2,2)}\ast s_{\l}$ is not multiplicity free for these $\l$'s.  If $\ell(\l)=2$
and $\l_1>\l_2>1$, then by Lemma~\ref{Lemma:2row} $s_{(n-2,2)}\ast s_{\l}$ is not
multiplicity free.  Thus, the only partitions left for consideration are
those in the statement of the theorem.

Since $s_{(n-2,2)}\ast s_{\l}$ is multiplicity free
if and only if $s_{(n-2,2)}\ast s_{\l'}$ is multiplicity free,  it  suffices to show
that $s_{(n-2,2)}\ast s_{\l}$ is multiplicity free when $\l=(n), (n-1,1)$, and
$\l=(m^k)$ where $m \geq k$. For $\l=(n)$ or $(1^n)$ the result is
trivial. For $\l=(n-1,1)$ or $(2,1^{n-2})$ the result follows from Proposition~\ref{P:p=1}. 

Let $m,k$ be positive integers  such that $n=mk$, $m \geq k$
and $k \geq 2$. Define the multiset
$$K_{(m^k),\alpha} := \{\nu\vdash n\, |\, k_{\a\n}^{(m^k)} \neq 0\}.$$
By \cite[Theorem~2.1]{vw}  or by the Skew Algorithm we have that $s_{(m^k)/(1,1)}
=s_{(m^{k-2},(m-1)^2)}$. Thus, by Pieri's rule we have
\begin{align*}
 K_{(m^2),(1,1)} &=  \{ (m,m), (m,m-1,1),(m-1,m-1,1,1)\}. \\
K_{(m^3),(1,1)}  &= \{(m+1,m,m-1),(m+1,m-1,m-1,1), (m^3),(m^2,m-1,1),
\\ 
&\kern8cm
(m,m-1,m-1, 1,1)\}.
\end{align*}
If $k>3$:
\begin{multline*}
K_{(m^k),(1,1)}=\{  ((m+1)^2,m^{k-4},(m-1)^2), (m+1, m^{k-2}, m-1),
 (m^k),  \\ 
 (m^{k-1}, m-1, 1),  (m+1, m^{k-3}, (m-1)^2, 1), (m^{k-2},(m-1)^2, 1^2)\}.
\end{multline*}

\begin{claim} Let $k \geq 2$. 
\medskip

\noindent (a) If
$m>3$, then 
$$K_{(m^k), (2)}=\{ (m+2, m^{k-2}, m-2), (m+1, m^{k-2}, m-2, 1),
(m^{k-1},m-2,2)\}.$$
\medskip

\noindent (b) $K_{(3^k),(2)} = \{(5,3^{k-2},1),
(4,3^{k-2},1^2)\}$.
\end{claim}

Assuming the claim, we have determined the diagrams that occur in each $K_{\a,(m^k)}$ when
$\a=(1,1)$ and $\a=(2)$ and they are disjoint, by Theorem~3.2, $s_{(n-2,2)}\ast
s_{(m^k)}$ is multiplicity free.
\end{proof}

\begin{proof}[Proof of Claim]
(a) From  \cite{d} we have that $g_{(n-2,2),(m^k),\n}=0$ if $\n_1
> m+2$. Hence $k_{(2), \n}^{(m^k)}=0$ and $k_{(1,1),\n}^{(m^k)}=0$ if $\n_1>m+2$.
Notice that if $\n_1< m$  the conditions of a Kronecker tableaux of shape
$(m^k)/(2)$ and type $\n/(2)$ cannot be satisfied. For example, let $\n_1=m-1$ and
consider the first two rows of a filling of $(m^k)/(2)$.  We must place two 2's in
the first row because there are not enough 1's to place two 1's in the second row.
This forces the following situation: 

\vskip10pt
\begin{picture}(400,25)(0,0)
\put(180,0){\line(1,0){90}} \put(200,20){\line(1,0){70}}
\put(180,10){\line(1,0){90}} \put(270,0){\line(0,0){20}}
\put(180,0){\line(0,0){10}}\put(190,0){\line(0,0){10}}\put(183,2){\scriptsize{1}}
\put(200,0){\line(0,0){20}}\put(210,0){\line(0,0){20}}\put(203,2){\scriptsize{2}}\put(218,2){$\cdots$}
\put(240,0){\line(0,0){20}}\put(250,0){\line(0,0){20}}\put(243,2){\scriptsize{2}}\put(253,2){\scriptsize{3}}
\put(263,2){\scriptsize{3}}\put(243,12){\scriptsize{1}}\put(253,12){\scriptsize{2}} \put(263,12){\scriptsize{2}}\put(260,0){\line(0,0){20}}\put(203,12){\scriptsize{1}}\put(218,12){$\cdots$}
\put(180,10){\red{\framebox(20,10)[t]{}}}\multiput(180,10)(.5,0){40}{\red{\line(0,1){10}}}
\end{picture}
\vskip10pt

In the empty box $(2,2)$ we cannot place another 1, and placing a 2
violates the $\alpha$-lattice permutation condition since there will be more 2's
than the number of 1's  plus 2 in the initial subword of length $2m-3$. A similar
argument can be applied if $\n_1<m-1$.  Therefore, if $k_{(2),\n}^{(m^k)}\in
K_{(m^k),(2)}$, we must have $m\leq \n_1\leq m+2$. \vskip 0inIf $\n_1 = m+2$, then
the first row of $(m^k)/(2)$ and the boxes $(2,1)$ and $(2,2)$ must be filled with
1's. This completely forces the following tableau: 

\vskip10pt
\begin{picture}(400,80)(0,0)
\put(180,10){\framebox(110,70)[t]{}}\put(180,20){\line(1,0){110}}\put(238,35){$\vdots$}
\put(180,30){\line(1,0){110}}\put(180,50){\line(1,0){110}}\put(195,10){\line(0,0){20}}
\put(180,60){\line(1,0){110}}\put(180,70){\line(1,0){110}}\put(195,50){\line(0,0){20}}
\put(210,10){\line(0,0){20}} \put(210,50){\line(0,0){30}}\put(225,10){\line(0,0){20}}
\put(225,50){\line(0,0){30}}\multiput(240,12)(0,10){2}{$\cdots$}\multiput(240,52)(0,10){3}{$\cdots$}
\put(275,10){\line(0,0){20}}\put(275,50){\line(0,0){30}}\multiput(180,70)(.5,0){60}{\red{\line(0,1){10}}}
\put(215,72){\scriptsize{1}}\put(280,72){\scriptsize{1}}\put(185,62){\scriptsize{1}} \put(200,62){\scriptsize{1}}
\put(215,62){\scriptsize{2}} \put(280,62){\scriptsize{2}} \put(185,52){\scriptsize{2}} \put(200,52){\scriptsize{2}}
\put(215,52){\scriptsize{3}} \put(280,52){\scriptsize{3}}\put(181,22){\scriptsize{$k$-2}}  \put(197,22){\scriptsize{$k$-2}}
 \put(212,22){\scriptsize{$k$-1}} \put(278,22){\scriptsize{$k$-1}} \put(181,12){\scriptsize{$k$-1}}
  \put(197,12){\scriptsize{$k$-1}} \put(215,12){\scriptsize{$k$}} \put(280,12){\scriptsize{$k$}}
\end{picture}
\vskip10pt

Hence, $k_{(2),(m+2,m^{k-2},m-2)}^{(m^k)}=1$. \vskip 0in Now suppose
that $\n_1=m+1$. In this case, we cannot place two 2's in the first row because
there will not be enough room for all the 1's.  Hence we must place two 1's in the
second row. This condition forces the following tableau:

\vskip10pt
\begin{picture}(400,80)(0,0)
\put(160,10){\framebox(140,70)[t]{}}\put(160,20){\line(1,0){140}}\put(238,35){$\vdots$}
\put(160,30){\line(1,0){140}}\put(160,50){\line(1,0){140}}\put(200,10){\line(0,0){20}}
\put(160,60){\line(1,0){140}}\put(160,70){\line(1,0){140}}\put(200,50){\line(0,0){30}}
\put(180,10){\line(0,0){20}} \put(180,50){\line(0,0){20}}\put(220,10){\line(0,0){20}}
\put(220,50){\line(0,0){30}}\multiput(235,12)(0,10){2}{$\cdots$}\multiput(235,52)(0,10){3}{$\cdots$}
\put(260,10){\line(0,0){20}}\put(260,50){\line(0,0){30}}\put(280,10){\line(0,0){20}}
\put(280,50){\line(0,0){30}}\multiput(160,70)(.5,0){80}{\red{\line(0,1){10}}}
\put(207,72){\scriptsize{1}}\put(267,72){\scriptsize{1}}\put(287,72){\scriptsize{2}}
\put(187,62){\scriptsize{1}}\put(167,62){\scriptsize{1}}\put(207,62){\scriptsize{2}}
\put(267,62){\scriptsize{2}}\put(287,62){\scriptsize{3}}\put(187,52){\scriptsize{2}}
\put(167,52){\scriptsize{2}}\put(207,52){\scriptsize{3}}\put(267,52){\scriptsize{3}}
\put(287,52){\scriptsize{4}}\put(163,22){\scriptsize{$k$-2}}\put(183,22){\scriptsize{$k$-2}}
\put(203,22){\scriptsize{$k$-1}}\put(263,22){\scriptsize{$k$-1}}\put(287,22){\scriptsize{$k$}}
\put(163,12){\scriptsize{$k$-1}}\put(183,12){\scriptsize{$k$-1}}
\put(207,12){\scriptsize{$k$}}\put(267,12){\scriptsize{$k$}}\put(283,12){\scriptsize{$k$+1}}
\end{picture}
\vskip10pt

Therefore, $k_{(2),(m+1, m^{k-2},m-2,1)}^{(m^k)}=1$. \vskip 0in Now if
$\n_1=m$, we cannot place all 1's in the first row of $(m^k)/(2)$, otherwise the
definition of Kronecker tableaux is not satisfied.  Hence we must place two 1's in
the second row and two 2's in the first row.  Now it is easy to see that in order
satisfy the condition of an $\alpha$-lattice permutation we are forced to fill the
rest of the diagram in the following unique way: 

\vskip10pt
\begin{picture}(400,80)(0,0)
\put(160,10){\framebox(160,70)[t]{}}\put(160,20){\line(1,0){160}}\put(233,33){$\vdots$}
\put(160,30){\line(1,0){160}}\put(160,50){\line(1,0){160}}\put(200,10){\line(0,0){20}}
\put(160,60){\line(1,0){160}}\put(160,70){\line(1,0){160}}\put(200,50){\line(0,0){30}}
\put(180,10){\line(0,0){20}} \put(180,50){\line(0,0){20}}\put(220,10){\line(0,0){20}}
\put(300,10){\line(0,0){20}}\put(300,50){\line(0,0){30}}
\put(220,50){\line(0,0){30}}\multiput(235,12)(0,10){2}{$\cdots$}\multiput(235,52)(0,10){3}{$\cdots$}
\put(260,10){\line(0,0){20}}\put(260,50){\line(0,0){30}}\put(280,10){\line(0,0){20}}
\put(280,50){\line(0,0){30}}\multiput(160,70)(.5,0){80}{\red{\line(0,1){10}}}
\put(207,72){\scriptsize{1}}\put(267,72){\scriptsize{1}}\put(287,72){\scriptsize{2}}
\put(307,72){\scriptsize{2}}\put(167,62){\scriptsize{1}}\put(187,62){\scriptsize{1}}
\put(207,62){\scriptsize{2}}\put(267,62){\scriptsize{2}}\put(287,62){\scriptsize{3}}
\put(307,62){\scriptsize{3}}\put(167,52){\scriptsize{2}}\put(187,52){\scriptsize{2}}
\put(207,52){\scriptsize{3}}\put(267,52){\scriptsize{3}}\put(287,52){\scriptsize{4}}
\put(307,52){\scriptsize{4}}\put(163,12){\scriptsize{$k$-2}}\put(183,12){\scriptsize{$k$-2}}
\put(203,22){\scriptsize{$k$-1}}\put(263,22){\scriptsize{$k$-1}}\put(287,22){\scriptsize{$k$}}
\put(307,22){\scriptsize{$k$}}\put(163,22){\scriptsize{$k$-1}}\put(183,22){\scriptsize{$k$-1}}
\put(207,12){\scriptsize{$k$}}\put(267,12){\scriptsize{$k$}}\put(283,12){\scriptsize{$k$+1}}
\put(303,12){\scriptsize{$k$+1}}
\end{picture}
\vskip10pt

Thus, $k_{(2),(m^{k-1},m-2,2)}^{(m^k)}=1$. 

\medskip
Part (b) is shown
in the same manner with the exception that it is not possible to obtain Kronecker
tableaux of shape $(m^k)/(2)$ and type $(m^{k-1},m-2,2)/(2)$ when $m=3$.
\end{proof}

\begin{remark} For  $n=4$, every partition is of the type listed in Theorem~\ref{T:p=2}. In
this case, although not covered by Theorem~\ref{T:p=2}, the result still holds that
$s_{(2,2)}\ast s_{\l}$ is multiplicity free for all $\l\vdash 4$. For
$n=5$, Theorem~4.5 
does not hold because in addition to the cases listed in the statement we also
have that $s_{(3,2)}\ast s_{(3,2)}$ and $s_{(3,2)}\ast s_{(2,2,1)}$ are multiplicity
free. 
\end{remark}

The following corollary follows directly from the proof of
Theorem~\ref{T:p=2}. Recall that $\chi(S)$ is the function that has value 1 if $S$
is true and 0 if $S$ is false.

\begin{cor} Let $n\geq 6$ and $n=mk$. Then
\begin{align*}
s_{(n-2,2)}\ast s_{(m^k)}&= s_{(m^k)}+ s_{(m^{k-1}, m-1, 1)}+
s_{(m^{k-2},(m-1)^2, 1^2)} + \chi(k\geq 4) s_{((m+1)^2,m^{k-4},(m-1)^2)}\\
& \kern1cm
+\chi(k\geq 3) s_{(m+1, m^{k-2}, m-1)}+\chi(k\geq 3) s_{(m+1, m^{k-3}, (m-1)^2, 1)} 
\\ 
&\kern1cm
+ s_{(m+2, m^{k-2}, m-2)}+ s_{(m+1, m^{k-2}, m-2, 1)}+ \chi(m\geq 4) s_{(m^{k-1},m-2,2)}.
\end{align*}
\end{cor}

\begin{lemma}
\label{Lemma:(p-1,1)}
Let $n$ and $p$ be positive integers such that $p\geq 3$ and let $\l\vdash n$ be such that
$\l_1\geq 2p-1$.  If $3\leq \ell(\l)$ or if $\ell(\l) =2$ and $\l_1>\l_2>1$, then
$$s_{(p-1,1)}s_{\l/(p-1,1)} -s_{(p-2,1)}s_{\l/(p-2,1)}$$
is not multiplicity free.
\end{lemma}
\begin{proof} Since $\l_1\geq 2p-1$, by Lemma~3.1.  the coefficient of $s_{\n}$ in the
Schur function expansion of $s_{(p- 1,1)}s_{\l/(p-1,1)} -s_{(p-2,1)}s_{\l/(p-2,1)}$
equals the number of Kronecker tableaux of shape \newline $\l/(p-1,1)$ and type
$\n/(p-1,1)$. To prove the lemma it suffices to show that it is possible to
construct two Kronecker tableaux of the same type.  We proceed by
cases. 

\medskip
\noindent {\sc Case 1.} Let $\ell(\l)\geq 3$ and $\l_2\geq p-1$.  Fill the rows of
$\l/(p-1,1)$ as follows:

{\leftskip1.5cm\parindent-15pt
\hangafter1
Row1:  Label box $(1,\l_1)$ with  2 and all other boxes with
  1's.
 
\hangafter1
Row 2:  Place  1's in boxes $(2, i)$ for $i=2, \ldots, p-1$. For $j=p, \ldots,
\l_2$, if box $(1,j)$ is labelled $x$ place $x+1$ in box $(2,j)$. 

\hangafter1
Row 3:
Place 1 in box $(3, 1)$. Then, if box $(2,i)$ is labelled $x$ place $x+1$ in box
$(3,i)$ for all $i=2, \ldots, \l_3$. 

\hangafter1
If $\ell(\l)\geq 4$, then place a 3 in
box $(4,1)$. For $i=2, \ldots,\l_4$,  place $x+1$ in box $(4, i)$ if box $(3,i)$
contains $x$. 

\hangafter1
For all other rows, place an $x+1$ under a box labelled $x$.
\par}

\medskip
Below we show a sample tableaux. In our example $\l_1=\l_2$ to illustrate the case
when we need to add a 3 in the second row.  If in addition $\l_2=\l_3$, then we
would also have a 4 at the end of the third row and so on.

\vskip10pt
\begin{picture}(400,50)
\put(150,40){\line(1,0){180}}\put(150,30){\line(1,0){180}}\put(150,20){\line(1,0){180}}
\put(150,10){\line(1,0){140}}\put(150,0){\line(1,0){110}}
\put(150,0){\line(0,0){40}}\multiput(150,30)(.5,0){120}{\red{\line(0,1){10}}}\put(210,0){\line(0,0){40}}
\multiput(150,20)(.5,0){20}{\red{\line(0,1){10}}}\put(160,00){\line(0,0){30}}
\put(220,0){\line(0,0){40}}\put(200,0){\line(0,0){30}}\put(170,0){\line(0,0){30}}
\put(250,0){\line(0,0){10}}\put(260,0){\line(0,0){10}}\put(280,10){\line(0,0){10}}
\put(290,10){\line(0,0){10}}\put(330,20){\line(0,0){20}}\put(320,20){\line(0,0){20}}
\put(310,20){\line(0,0){20}}
\put(153,12){\scriptsize{1}}\put(153,2){\scriptsize{3}}
\put(163,22){\scriptsize{1}}\put(163,12){\scriptsize{2}}\put(163,2){\scriptsize{3}}
\multiput(178,2)(0,10){3}{$\cdots$}\put(203,22){\scriptsize{1}}\put(203,12){\scriptsize{2}}
\put(203,2){\scriptsize{3}}\put(213,32){\scriptsize{1}}\put(213,22){\scriptsize{2}}
\put(213,12){\scriptsize{3}}\put(213,2){\scriptsize{4}}\put(228,2){$\cdots$}
\put(240,12){$\cdots$}\put(260,22){$\cdots$}\put(260,32){$\cdots$}
\put(253,2){\scriptsize{4}}\put(283,12){\scriptsize{3}}\put(313,32){\scriptsize{1}}
\put(323,32){\scriptsize{2}}\put(313,22){\scriptsize{2}}\put(323,22){\scriptsize{3}}
\put(120,15){$T=$}
\end{picture}
\vskip10pt

The tableau $T$ is a SSYT with $p-2$ boxes labelled $1$ in the second
row. It is straightforward to check that its reverse reading word is a
$(p-1,1)$-lattice permutation.  Hence, $T$ is a Kronecker tableau of shape
$\l/(p-1,1)$.

We can switch the 2 in box $(1,\l_1)$ with the 1 in box $(3,1)$ to obtain
a different  Kronecker tableau of the same shape and type as $T$. Hence,
$k_{(p-1,1),type(T)}^\l\geq 2$. 

\medskip
\noindent {\sc Case 2.} Let $\ell(\l)\geq 3$
and $2\leq\l_2<p-1$. In this case fill the diagram $\l/(p-1,1)$ as follows:

{\leftskip1.5cm\parindent-15pt
\hangafter1
Row 1: Place 2's in the last $p-2$ boxes  and  1's in all other boxes. Since
$\l_1\geq 2p-1$, we always have at least two 1's in the first row. 

\hangafter1
Row 2:
Place 1's in boxes $(2,i)$ for $i=2, \ldots, \l_2-1$ (if $\l_2>2$) and place 2 in
box $(2,\l_2)$. (If $\l_2=2$, there are no 1's in the second row.) Notice that
placing the 2 in this box does not violate the definition of Kronecker tableaux
because in the initial subword of length $\l_1-p+2$ we have $\# 1's +(p-1)\geq p+1$
and  $\#2's +1 =p$.

\hangafter1
Row 3: Place a 1 in box $(3,1)$  and 3's in the remaining boxes. 

\hangafter1
If
$\ell(\l)\geq 4$, then in row 4 place a 3 in box $(4,1)$ and 4's in the remaining
boxes. 

\hangafter1
For all remaining rows place an $x+1$ under a box labelled $x$.
\par}

\medskip
We illustrate the first four rows of such a tableau:

\vskip10pt
\begin{picture}(400,50)
\put(150,40){\line(1,0){220}}\put(150,30){\line(1,0){220}}\put(150,20){\line(1,0){120}}
\put(150,10){\line(1,0){80}}\put(150,0){\line(1,0){50}}\put(150,0){\line(0,0){40}}
\put(160,0){\line(0,0){30}}\multiput(150,30)(.5,0){280}{\red{\line(0,1){10}}}
\multiput(150,20)(.5,0){20}{\red{\line(0,1){10}}}\put(170,0){\line(0,0){30}}
\put(190,0){\line(0,0){10}}\put(200,0){\line(0,0){10}}\put(230,10){\line(0,0){10}}
\put(220,10){\line(0,0){10}}\put(270,20){\line(0,0){10}}\put(260,20){\line(0,0){10}}
\put(250,20){\line(0,0){10}}\put(290,30){\line(0,0){10}}\put(300,30){\line(0,0){10}}
\put(320,30){\line(0,0){10}}\put(330,30){\line(0,0){10}}\put(340,30){\line(0,0){10}}
\put(360,30){\line(0,0){10}}\put(370,30){\line(0,0){10}}
\put(153,2){\scriptsize{3}}\put(153,12){\scriptsize{1}}
\put(163,2){\scriptsize{4}}\put(163,12){\scriptsize{3}}\put(163,22){\scriptsize{1}}
\put(175,2){$\cdots$}\put(190,12){$\cdots$}\put(200,22){$\cdots$}\put(305,32){$\cdots$}
\put(193,2){\scriptsize{4}}\put(223,12){\scriptsize{3}}
\put(253,22){\scriptsize{1}}\put(263,22){\scriptsize{2}}
\put(293,32){\scriptsize{1}}\put(323,32){\scriptsize{1}}
\put(333,32){\scriptsize{2}}\put(343,32){$\cdots$}
\put(363,32){\scriptsize{2}}\put(120,20){$T'=$}
\end{picture}
\vskip10pt

The tableau $T'$ is a Kronecker tableaux of shape $\l/(p-1,1)$.
Switching the 2 in box $(2,\l_2)$ with the 1 in box $(3,1)$ yields another Kronecker
tableaux of the same type and shape as $T'$. Hence, $k_{(p-1,1),type(T')}^{\l}\geq
2$.

\medskip
\noindent {\sc Case 3.} Let $\ell(\l)\geq 3$ and $\l_2=1$. In this case
$\l$ is a hook partition. We illustrate two Kronecker tableaux of shape $\l/(p-1,1)$
and  type $(\l_1-p+2, p-1, 1^{n-\l_1-1})/(p-1,1)$: 

\vskip10pt
\centerline{\begin{picture}(400,80)
\put(10,80){\line(1,0){120}}\put(10,70){\line(1,0){120}}
\multiput(10,70)(.5,0){80}{\red{\line(0,1){10}}}\multiput(10,60)(.5,0){20}{\red{\line(0,1){10}}}
\put(10,10){\line(0,0){70}}\put(20,10){\line(0,0){60}}\put(10,60){\line(1,0){10}}
\put(10,50){\line(1,0){10}}\put(10,40){\line(1,0){10}}\put(14,24){$\vdots$}
\put(10,10){\line(1,0){10}}\put(10,20){\line(1,0){10}}
\put(50,70){\line(0,0){10}}\put(60,70){\line(0,0){10}}
\put(80,70){\line(0,0){10}}\put(90,70){\line(0,0){10}}
\put(100,70){\line(0,0){10}}\put(120,70){\line(0,0){10}}
\put(130,70){\line(0,0){10}}
\put(13,52){\scriptsize{3}}\put(13,42){\scriptsize{4}}\put(13,12){\scriptsize$k$}
\put(53,72){\scriptsize{1}}\put(65,72){$\cdots$}\put(83,72){\scriptsize{1}}
\put(93,72){\scriptsize{2}}\put(105,72){$\cdots$}\put(123,72){\scriptsize{2}}
\put(210,80){\line(1,0){130}}\put(210,70){\line(1,0){130}}
\multiput(210,70)(.5,0){80}{\red{\line(0,1){10}}}\multiput(210,60)(.5,0){20}{\red{\line(0,1){10}}}
\put(210,10){\line(0,0){70}}\put(220,10){\line(0,0){60}}\put(210,60){\line(1,0){10}}
\put(210,50){\line(1,0){10}}\put(210,40){\line(1,0){10}}\put(214,24){$\vdots$}
\put(210,10){\line(1,0){10}}\put(210,20){\line(1,0){10}}
\put(250,70){\line(0,0){10}}\put(260,70){\line(0,0){10}}
\put(280,70){\line(0,0){10}}\put(290,70){\line(0,0){10}}
\put(300,70){\line(0,0){10}}\put(320,70){\line(0,0){10}}
\put(330,70){\line(0,0){10}} \put(340,70){\line(0,0){10}}
\put(213,52){\scriptsize{1}}\put(213,42){\scriptsize{4}}\put(213,12){\scriptsize$k$}
\put(253,72){\scriptsize{1}}\put(265,72){$\cdots$}\put(283,72){\scriptsize{1}}
\put(293,72){\scriptsize{2}}\put(305,72){$\cdots$}\put(323,72){\scriptsize{2}}
\put(333,72){\scriptsize{3}}
\end{picture}}
\vskip10pt

Here $\ell(\l)=k$ and  if there are $m$ 1's in the tableau on the left hand side,
then there are $m-1$ 1's in  the first row of the tableau on the right. The number
of 2's in the first row is $p-2$. 

By cases (1)-(3) we have that
$s_{(p-1,1)}\ast s_{\l/(p-1,1)} -s_{(p-2,1)}\ast s_{\l/(p-2,1)}$ is not multiplicity
free if $\ell(\l)\geq 3$.

\medskip
\noindent {\sc Case 4.}  Let $\ell(\l)=2$ and $\l_1\neq \l_2$ and $\l_2>1$.  We have two
subcases. 

\medskip
\noindent Case (i): $\l_2 > p-1$.  In this case the following
is a Kronecker tableau of shape $\l/(p-1,1)$ and type $\n/(p-1,1)$, where
$\n=(\l_1+(p-3),\l_2-(p-2),1)$: 

\vskip10pt
\begin{picture}(400,30)
\put(150,30){\line(1,0){150}}
\put(150,20){\line(1,0){150}}
\put(150,10){\line(1,0){110}}
\put(150,10){\line(0,0){20}}
\multiput(150,20)(.5,0){120}{\red{\line(0,1){10}}}\multiput(150,10)(.5,0){20}{\red{\line(0,1){10}}}
\put(160,10){\line(0,0){10}}\put(170,10){\line(0,0){10}}\put(180,12){$\cdots$}
\put(200,10){\line(0,0){10}}\put(210,10){\line(0,0){20}}\put(220,10){\line(0,0){20}}
\put(230,12){$\cdots$}\put(250,10){\line(0,0){10}}\put(260,10){\line(0,0){10}}
\put(250,22){$\cdots$}\put(280,20){\line(0,0){10}}\put(290,20){\line(0,0){10}}
\put(300,20){\line(0,0){10}}
\put(163,12){\scriptsize{1}}\put(203,12){\scriptsize{1}}\put(213,12){\scriptsize{2}}
\put(213,22){\scriptsize{1}}\put(253,12){\scriptsize{2}}\put(283,22){\scriptsize{1}}
\put(293,22){\scriptsize{3}}
\end{picture}
\vskip10pt

We obtain another Kronecker tableau of the same shape and type by
switching the 3 in box $(1,\l_1)$ with the 2 in box $(2,\l_2)$. Since $\l_1\neq
\l_2$, this is always possible. Hence,  $k_{(p-1,1),\n}^\l\geq 2$. 

\medskip
\noindent Case (ii): $2\leq \l_2\leq p-1$. In this case the following are Kronecker
tableaux of shape $\l/(p-1,1)$ and type $\n/(p-1,1)$, where
$\n=(\l_1+\l_2-p,p-1,1)$: 

\vskip10pt
\begin{picture}(400,30)
\put(50,20){\line(1,0){160}}
\put(50,10){\line(1,0){160}}
\put(50,0){\line(1,0){60}}
\put(50,0){\line(0,0){20}}\put(60,0){\line(0,0){10}}
\multiput(50,10)(.5,0){160}{\red{\line(0,1){10}}}\multiput(50,0)(.5,0){20}{\red{\line(0,1){10}}}
\put(70,0){\line(0,0){10}}\put(90,0){\line(0,0){10}}\put(100,0){\line(0,0){10}}
\put(110,0){\line(0,0){10}}\put(130,10){\line(0,0){10}}\put(140,10){\line(0,0){10}}
\put(160,10){\line(0,0){10}}\put(170,10){\line(0,0){10}}\put(180,10){\line(0,0){10}}
\put(200,10){\line(0,0){10}}\put(210,10){\line(0,0){10}}
\put(63,2){\scriptsize{1}}\put(75,2){$\cdots$}\put(93,2){\scriptsize{1}}\put(103,2){\scriptsize{3}}
\put(133,12){\scriptsize{1}}\put(145,12){$\cdots$}\put(163,12){\scriptsize{1}}
\put(173,12){\scriptsize{2}}\put(185,12){$\cdots$}\put(203,12){\scriptsize{2}}
\put(250,20){\line(1,0){170}}
\put(250,10){\line(1,0){170}}
\put(250,0){\line(1,0){60}}
\put(250,0){\line(0,0){20}}\put(260,0){\line(0,0){10}}
\multiput(250,10)(.5,0){160}{\red{\line(0,1){10}}}\multiput(250,0)(.5,0){20}{\red{\line(0,1){10}}}
\put(270,0){\line(0,0){10}}\put(290,0){\line(0,0){10}}\put(300,0){\line(0,0){10}}
\put(310,0){\line(0,0){10}}\put(330,10){\line(0,0){10}}\put(340,10){\line(0,0){10}}
\put(360,10){\line(0,0){10}}\put(370,10){\line(0,0){10}}\put(380,10){\line(0,0){10}}
\put(400,10){\line(0,0){10}}\put(410,10){\line(0,0){10}}\put(420,10){\line(0,0){10}}
\put(263,2){\scriptsize{1}}\put(275,2){$\cdots$}\put(293,2){\scriptsize{1}}\put(303,2){\scriptsize{1}}
\put(333,12){\scriptsize{1}}\put(345,12){$\cdots$}\put(363,12){\scriptsize{1}}
\put(373,12){\scriptsize{2}}\put(385,12){$\cdots$}\put(403,12){\scriptsize{2}}
\put(413,12){\scriptsize{3}}
\end{picture}
\vskip10pt

There are $p-2$ boxes labelled 2 in the first row of these diagrams. In
the first row, in the left tableau there is one more 1 than in the right tableau.
Hence, $k_{(p-1,1),\n}^{(\l_1,\l_2)}\geq 2$ whenever $\l_1>\l_2>1$.  This concludes
the proof.
 \end{proof}

 \begin{theorem}
 Let $n>16$ and $\l\vdash n$. Then $s_{(n-3,3)}\ast s_{\l}$ is multiplicity free if and only if
 $\l=(n), (1^n), (n-1,1),(2,1^{n-2})$ and if $n$ is even also $\l=(n/2, n/2)$ or $(2^{n/2})$.
 \end{theorem}
 \begin{proof}
 Since $n>16$, we  have  $\ell(\l)\geq 5$ or $\l_1\geq 5$.
 Without loss of generality, for the remainder of the proof we can assume that $\l_1\geq 5$.
  The results
 for $(1^n)$, $(2, 1^{n-2})$ and $(2^{n/2})$ will follow by the
 fact that they are conjugate to $(n), (n-1,1)$
 and $(n/2, n/2)$ respectively.

Since $\l_1\geq 5$,  by Theorem~3.2 we have
$g_{(n-3,3),\n}^{\l}=k_{(3),\n}^{\l}+k_{(2,1),\n}^\l+k_{(1,1,1),\n}^\l$.
 By Lemma~\ref{Lemma:free}, if $\ell(\l)>3$  and $\l$ is not a rectangular
 partition, then there exists a $\n$ such that $k_{(1,1,1),\n}^{\l}\geq 2$.
 Also,  by Lemma~\ref{Lemma:(p-1,1)}, if $\ell(\l)\geq 3$ or if $\ell(\l)=2$ and
 $\l_1>\l_2>1$,
  then there exists a $\n$ such
 that $k_{(2,1),\n}^\l\geq 2$. Hence, for all these cases $s_{(n-3,3)}\ast s_{\l}$
 is not multiplicity free.
 It remains to  investigate the  cases in which
 $\ell(\l)\leq 2$ and either $\l_2=0,1$ or $n/2$ (when $n$ is even).
 If $\l_2=0$, it is trivially true that $s_{(n-3,3)}\ast s_{(n)}$ is multiplicity free.
 If $\l_2=1$, since $C(n-3,3)=1$, it follows by
 Proposition~\ref{P:p=1} that $s_{(n-3,3)}\ast s_{(n-1,1)}$ is multiplicity free.

 Now assume that $n$ is even and $\l=(n/2,n/2)$.  Recall that $K_{\a,\l}$ is the multiset
 containing all  $\n$'s for which  $k_{\a\n}^\l\neq 0$. The multiplicity of each $\n$
 in $K_{\a,\l}$ is $k_{\a\n}^\l$. We have that $K_{(1,1,1),(n/2,n/2)}=\emptyset$,
 \begin{multline*}
 K_{(2,1),(\frac{n}{2},\frac{n}{2})} = \left\{ \left(\frac{n}{2}+1, \frac{n}{2}-1\right),
 \left(\frac{n}{2}+1, \frac{n}{2}-2, 1\right), \left(\frac{n}{2}, \frac{n}{2}-1,1\right),
 \left(\frac{n}{2},\frac{n}{2}-2,1,1\right), \right. \\\left. 
 \left(\frac{n}{2},\frac{n}{2}-2,2\right),\left(\frac{n}{2}-1, \frac{n}{2}-1, 2\right), \left(\frac{n}{2}-1, \frac{n}{2}-2, 2, 1\right)\right\},
 \end{multline*}
 and
\begin{multline*} 
K_{(3),(n/2,n/2)}=\left\{\left(\frac{n}{2}+3, \frac{n}{2}-3\right), \left(\frac{n}{2}+2, \frac{n}{2}-3,1\right), \right.\\\left.
\left(\frac{n}{2}+1,
 \frac{n}{2}-3, 2\right),\left(\frac{n}{2}, \frac{n}{2}-3,3\right)\right\}.
\end{multline*}

 By inspection we see that
 $K_{(2,1),(\frac{n}{2},\frac{n}{2})}$ and $K_{(3),(\frac{n}{2},\frac{n}{2})}$ are disjoint.  Furthermore,
 we have
\begin{multline*}
 s_{(n-3,3)}\ast s_{(\frac{n}{2},\frac{n}{2})} = s_{(\frac{n}{2}+1, \frac{n}{2}-1)}+
 s_{(\frac{n}{2}+1, \frac{n}{2}-2, 1)}+ s_{(\frac{n}{2}, \frac{n}{2}-1,1)}+
 s_{(\frac{n}{2}, \frac{n}{2}-2,1,1)}
 +s_{(\frac{n}{2},\frac{n}{2}-2,2)}+s_{(\frac{n}{2}, \frac{n}{2}-3,3)}\\
 +s_{(\frac{n}{2}-1, \frac{n}{2}-1, 2)}+
 s_{(\frac{n}{2}-1, \frac{n}{2}-2, 2, 1)}+s_{(\frac{n}{2}+3, \frac{n}{2}-3)}+
  s_{(\frac{n}{2}+2, \frac{n}{2}-3,1)}+ s_{(\frac{n}{2}+1, \frac{n}{2}-3, 2)}.
\end{multline*}
 \end{proof}
 \begin{theorem}
 Let $p\geq 4$ and $n> (2p-2)^2$. Then $s_{(n-p,p)}\ast s_{\l}$ is multiplicity free if and only if $\l=(n), (1^n), (n-1,1),$ or $(2, 1^{n-2})$.
 \end{theorem}
 \begin{proof}  Since $n>(2p-2)^2$ we  have  $\l_1$ or $\ell(\l)$  greater or equal to $2p-1$.
 Without loss of generality we  assume that $\l_1\geq 2p-1$.
 It is well known that if $\l=(n)$ or $(1^n)$,
 then $s_{(n-p ,p)}\ast s_{\l}$ is multiplicity free. Since $C(n-p,p)\leq 1$,
 by Proposition~\ref{P:p=1} and Corollary~\ref{cor:cor2},
 $s_{(n-p,p)}\ast s_{\l}$ is also multiplicity free when $\l=(n-1,1)$ or $(2, 1^{n-2})$.

  By Lemma~\ref{Lemma:free} and Lemma~\ref{Lemma:(p-1,1)},  if $\l_1\geq 2p-1$, then
   $s_{(n-p,p)}\ast s_{\l}$ could be multiplicity free only when $\l_2 = 0, 1$ or,
    if $n$ is even,
 when $\l=(\frac{n}{2},\frac{n}{2})$.  Hence, we only have to show that for
 $p\geq4$, $s_{(n-p,p)}\ast s_{(\frac{n}{2}, \frac{n}{2})}$ is not multiplicity free.
 We consider  the case  $p=4$ separately. The following are two Kronecker tableaux
 of shape $(n/2,n/2)/\a$ and type $\n/\a$
 where $\n= (\frac{n}{2}, \frac{n}{2}-2, 2)$ and $\a=(3,1)$ and $\a=(2,2)$ respectively.

 \vskip10pt
 \begin{picture}(400,30)
 \put(50,0){\framebox(100,20)[t]{}}\put(50,10){\line(1,0){100}}
 \multiput(50,10)(.5,0){60}{\red{\line(0,1){10}}}\multiput(50,0)(.5,0){20}{\red{\line(0,1){10}}}
 \put(60,0){\line(0,0){20}}\put(70,0){\line(0,0){20}}\put(80,0){\line(0,0){20}}
 \put(90,0){\line(0,0){20}}\put(120,0){\line(0,0){20}}\put(130,0){\line(0,0){20}}
 \put(140,0){\line(0,0){20}}
 \put(63,2){\scriptsize{1}}\put(73,2){\scriptsize{1}}\put(83,2){\scriptsize{2}}\put(100,2){$\cdots$}
 \put(123,2){\scriptsize{2}}\put(133,2){\scriptsize{3}} \put(143,2){\scriptsize{3}}
 \put(83,12){\scriptsize{1}} \put(100,12){$\cdots$}\put(123,12){\scriptsize{1}}
 \put(133,12){\scriptsize{2}}\put(143,12){\scriptsize{2}}
 \put(200,0){\framebox(100,20)[t]{}}\put(200,10){\line(1,0){100}}
 \multiput(200,10)(.5,0){40}{\red{\line(0,1){10}}}\multiput(200,0)(.5,0){40}{\red{\line(0,1){10}}}
 \put(210,0){\line(0,0){20}}\put(220,0){\line(0,0){20}}\put(230,0){\line(0,0){20}}
 \put(240,0){\line(0,0){20}}\put(270,0){\line(0,0){20}}\put(280,0){\line(0,0){20}}
 \put(290,0){\line(0,0){20}}
\put(223,2){\scriptsize{2}}\put(233,2){\scriptsize{2}}\put(250,2){$\cdots$}\put(273,2){\scriptsize{2}}
 \put(283,2){\scriptsize{3}} \put(293,2){\scriptsize{3}}\put(223,12){\scriptsize{1}}
 \put(233,12){\scriptsize{1}} \put(250,12){$\cdots$}\put(273,12){\scriptsize{1}}
 \put(283,12){\scriptsize{1}}\put(293,12){\scriptsize{1}}
  \end{picture}
 \vskip10pt

 Hence, $k_{(3,1),\n}^{(\frac{n}{2},\frac{n}{2})} +k_{(2,2),\n}^{(\frac{n}{2},\frac{n}{2})}\geq 2$. Therefore,
 by Theorem~3.2, $g_{(n-4,4),(\frac{n}{2},\frac{n}{2}),\n}\geq 2$.

\smallskip
 If $p\geq 5$, the following are Kronecker tableaux of shape $(n/2,n/2)/\a$ and type $\n/\a$,
where
 $\n=(\frac{n}{2}+p-4, \frac{n}{2}-p+2, 2)$ and $\a=(p-1,1)$ and $\a=(p-2,2)$ respectively:

 \vskip10pt
 \begin{picture}(400,30)
\put(50,0){\framebox(130,20)[t]{}}\put(50,10){\line(1,0){130}}
\multiput(50,10)(.5,0){120}{\red{\line(0,1){10}}}\multiput(50,0)(.5,0){20}{\red{\line(0,1){10}}}
\put(60,0){\line(0,0){10}}\put(70,0){\line(0,0){10}}\put(100,0){\line(0,0){10}}
\put(110,0){\line(0,0){20}}\put(120,0){\line(0,0){20}}\put(150,0){\line(0,0){20}}
\put(160,0){\line(0,0){20}}\put(170,0){\line(0,0){20}}
\put(63,2){\scriptsize{1}}\put(80,2){$\cdots$}\put(103,2){\scriptsize{1}}\put(113,2){\scriptsize{2}}
\put(130,2){$\cdots$}\put(153,2){\scriptsize{2}}\put(163,2){\scriptsize{3}}\put(173,2){\scriptsize{3}}
\put(113,12){\scriptsize{1}}\put(130,12){$\cdots$}\put(153,12){\scriptsize{1}}\put(163,12){\scriptsize{2}}
\put(173,12){\scriptsize{2}}
\put(250,0){\framebox(130,20)[t]{}}\put(250,10){\line(1,0){130}}
\multiput(250,10)(.5,0){120}{\red{\line(0,1){10}}}\multiput(250,0)(.5,0){40}{\red{\line(0,1){10}}}
\put(260,0){\line(0,0){10}}\put(270,0){\line(0,0){10}}\put(280,0){\line(0,0){10}}
\put(300,0){\line(0,0){10}}
\put(310,0){\line(0,0){20}}\put(320,0){\line(0,0){20}}\put(350,0){\line(0,0){20}}
\put(360,0){\line(0,0){20}}\put(370,0){\line(0,0){20}}
\put(273,2){\scriptsize{1}}\put(285,2){$\cdots$}\put(303,2){\scriptsize{1}}\put(313,2){\scriptsize{2}}
\put(330,2){$\cdots$}\put(353,2){\scriptsize{2}}\put(363,2){\scriptsize{3}}\put(373,2){\scriptsize{3}}
\put(313,12){\scriptsize{1}}\put(330,12){$\cdots$}\put(353,12){\scriptsize{1}}\put(363,12){\scriptsize{1}}
\put(373,12){\scriptsize{2}}
 \end{picture}
 \vskip10pt

 Hence, we have $k_{(p-1,1),\n}^{(\frac{n}{2},\frac{n}{2})} +k_{(p-2,2),\n}^{(\frac{n}{2},\frac{n}{2})}\geq 2$.
 Therefore, by Theorem~3.2, $g_{(n-p,p),(\frac{n}{2},\frac{n}{2}),\n}\geq 2$ for all $p\geq 5$.
 \end{proof}


\subsection{The Kronecker Product of a two row shape and a hook shape}

In this subsection we show how to use the combinatorial rule of Theorem~3.2  to
obtain formulas for the coefficients in the Kronecker product $s_{(n-p,p)}*s_{(n-s,
1^s)}$, where $n-s \geq 2p-1$. The formulas obtained are equivalent to the results
obtained by Remmel \cite{r1}. For this reason, we include only one example: the
multiplicity of a hook in $s_{(n-p,p)}*s_{(n-s, 1^s)}$ in the case $p\geq 2$ and
$n-s \geq 2p-1$. 

We first need the following preliminary result.

\begin{prop}
\label{P:bds} Let $s,p,n$ be positive integers such that $n-s\geq 2p-1$ and let
$\n=(\n_1,\n_2,\ldots, \n_{\ell(\n)})\vdash n$.
\begin{enumerate}
\item[(a)] If $g_{(n-p,p),(n-s,1^s),\n} \neq 0$, then $n-s-p\leq \n_1\leq n-s+1$, $\max\{1,p-s\}\leq \n_2\leq p+1$,  and $\n_i \leq 2$ for $i\geq 3$.
\item[(b)] $g_{(n-p,p),(n-s,1^s),\n} = k_{(\n_2,1^{p-\n_2}),\n}^{(n-s,1^s)}+k_{(\n_2-1,1^{p-\n_2+1}),\n}^{(n-s,1^s)}$.
\end{enumerate}
\end{prop}

\begin{proof} (a) Given $\n=(\n_1, \n_2, \ldots, \n_{\ell(\n)}) \vdash n$, let $\a \vdash
p$ such that $\a \subseteq (n-s,1^s)$ and $\a \subseteq \n$. Thus $\a=(p-u,1^u)$, $u
\leq \min\{s, \ell(\n)-1\}$, $p-u \leq \n_1$. The maximum number of $1$'s in a
Kronecker tableau occurs when $\a=(1^p)$ (if $p \leq \min\{s, \ell(\n)-1\}$) and the
first row is filled with $1$'s and there is a $1$ in the first column of $\l/\a$.
Therefore $\n_1 \leq n-s+1$. The minimum number of $1$'s occurs when $\a=(p)$  and
there are no $1$'s in the first column of $\l/\a$. Therefore $\n_1 \geq n-s-p$. The
maximum number of $2$'s occurs when $\a=(p)$ and there is one $2$ in the first
column of $\l/\a$. Thus $\n_2 \leq p+1$. If $p \leq s+1$, the minimum number of
$2$'s occurs when $\a=(1^p)$ and it equals 1. In this case there are no $2$'s in
$\l/\a$. If $p
>s+1$, the minimum number of $2$'s occurs when $\a=(p-s,s)$ and it equals $p-s$ (no
$2$'s in the first column of $\l/\a$). Thus $\n_2 \geq \max\{1, p-s\}$.  If $i \geq
3$, $\a_i \leq 1$. If $\a_i=1$ at most one $i$ can be placed in the first column of
$\l/\a$. If $\a_i=0$ at most one $i$ can be placed in the fist row of $\l/\a$ and at
most one $i$ can be placed in the first column of $\l/\a$. Thus $\n_i \leq 2$ for $i
\geq 3$.

(b) When forming Kronecker tableaux of shape $\l/\a$ and type $\n/\a$ we cannot fill
the box in position $(2,\a_1)$ with $1$. If $\a_1 \neq 1$ or $p$, we must place
$2$'s in exactly  $\a_1-1$ boxes of the first row of $\l/\a$ and we can place at
most one $2$ in the first column of $\l/\a$. Since $\a_2=1$, we have $\n_2=\a_1$ or
$\n_2=\a_1+1$. If $\a_1= p$, we must place $2$'s in exactly  $\a_1=p$ boxes of the
first row of $\l/\a$ and we can place at most one $2$ in the first column of
$\l/\a$. Since $\a_2=0$, we have  $\n_2=\a_1$ or $\n_2=\a_1+1$. If $\a_1 =1$, we can
place at most one $2$ in the first column of $\l/\a$. Since $\a_2=1$, we have
$\n_2=1= \a_1$ or $\n_2=2= \a_1+1$.
\end{proof}

\begin{note} If $n-s=2p-1$, then $k_{(p), \n}^{(n-s,1^s)}=0$ and if $s_\n$
appears in the decomposition of $s_{(n-p,p)}*s_{(n-s,1^s)}$ then  $\n_2 \leq p$ and
$\n_1 \geq n-s-p+1=p$. 
\end{note}

By analyzing Kronecker tableaux we have obtained
the following well-known result \cite{r1}.

\begin{cor}
\label{Cor:hkdh} If $g_{(n-p,p),(n-s,1^s),\n}\neq 0$, then $\n$ is a hook,
$\n=(\n_1,1^{n-\n_1})$, or a double hook, $\n=(\n_1,\n_2,2^i,1^j)$.

\end{cor}

To show the applicability of Theorem~3.2, we use the previous proposition to obtain
the coefficient of a hook in the Kronecker product $s_{(n-p,p)}*s_{(n-s, 1^s)}$.
\medskip

It follows from Proposition~\ref{P:bds} that
$g_{(n-p,p),(n-s,1^s),(n-t,1^t)}=k_{(1^p), (n-t,1^t)}^{(n-s,1^s)}$.
\medskip

To obtain the possible Kronecker tableaux  we must have $p \leq s+1$ and no boxes of
$(n-s,1^s)/(1^p)$ can be filled with $2$'s. 

\vskip10pt
\begin{picture}(400,80)
\put(30,70){\framebox(60,10)[t]{}}\put(45,70){\line(0,0){10}}\put(35,72){\scriptsize{1}}
\put(53,72){$\cdots$}\put(75,70){\line(0,0){10}}\put(80,72){\scriptsize{1}}
\put(15,0){\line(0,0){80}}\put(30,0){\line(0,0){70}}\put(15,80){\line(1,0){15}}
\multiput(15,50)(0,.5){60}{\red{\line(1,0){15}}}\put(15,50){\line(1,0){15}}
\put(20,42){\scriptsize{1}}\put(15,40){\line(1,0){15}}\put(15,32){\scriptsize{$p$+1}}
\put(15,30){\line(1,0){15}}\put(20,15){$\vdots$}\put(15,10){\line(1,0){15}}
\put(15,0){\line(1,0){15}}\put(20,2){\scriptsize{$s$}}
\put(-5,-15){$\n=(n-s+1,1^{s-1})$}
\put(150,70){\framebox(60,10)[t]{}}\put(165,70){\line(0,0){10}}\put(155,72){\scriptsize{1}}
\put(173,72){$\cdots$}\put(195,70){\line(0,0){10}}\put(200,72){\scriptsize{1}}
\put(135,0){\line(0,0){80}}\put(150,0){\line(0,0){70}}\put(135,80){\line(1,0){15}}
\multiput(135,50)(0,.5){60}{\red{\line(1,0){15}}}\put(135,50){\line(1,0){15}}
\put(135,42){\scriptsize{$p$+1}}\put(135,40){\line(1,0){15}}\put(135,32){\scriptsize{$p$+2}}
\put(135,30){\line(1,0){15}}\put(141,15){$\vdots$}\put(135,10){\line(1,0){15}}
\put(135,0){\line(1,0){15}}\put(135,2){\scriptsize{$s$+1}}
\put(120,-15){$\n=(n-s,1^{s})$}
\put(270,70){\framebox(70,10)[t]{}}\put(280,70){\line(0,0){10}}\put(272,72){\scriptsize{1}}
\put(290,72){$\cdots$}\put(310,70){\line(0,0){10}}\put(314,72){\scriptsize{1}}
\put(323,70){\line(0,0){10}}\put(325,72){\scriptsize{$p$+1}}
\put(255,0){\line(0,0){80}}\put(270,0){\line(0,0){70}}\put(255,80){\line(1,0){15}}
\multiput(255,50)(0,.5){60}{\red{\line(1,0){15}}}\put(255,50){\line(1,0){15}}
\put(260,42){\scriptsize{1}}\put(255,40){\line(1,0){15}}\put(255,32){\scriptsize{$p$+2}}
\put(255,30){\line(1,0){15}}\put(261,15){$\vdots$}\put(255,10){\line(1,0){15}}
\put(255,0){\line(1,0){15}}\put(255,2){\scriptsize{$s$+1}}
\put(240,-15){$\n=(n-s,1^{s})$}
\put(370,70){\framebox(70,10)[t]{}}\put(380,70){\line(0,0){10}}\put(372,72){\scriptsize{1}}
\put(390,72){$\cdots$}\put(410,70){\line(0,0){10}}\put(414,72){\scriptsize{1}}
\put(423,70){\line(0,0){10}}\put(425,72){\scriptsize{$p$+1}}
\put(355,0){\line(0,0){80}}\put(370,0){\line(0,0){70}}\put(355,80){\line(1,0){15}}
\multiput(355,50)(0,.5){60}{\red{\line(1,0){15}}}\put(355,50){\line(1,0){15}}
\put(355,42){\scriptsize{$p$+2}}\put(355,40){\line(1,0){15}}\put(355,32){\scriptsize{$p$+3}}
\put(355,30){\line(1,0){15}}\put(361,15){$\vdots$}\put(355,10){\line(1,0){15}}
\put(355,0){\line(1,0){15}}\put(355,2){\scriptsize{$s$+1}}
\put(340,-15){$\n=(n-s-1,1^{s+1})$}
\end{picture}
\vskip20pt

Hence, the highest possible coefficient in this case is 2 and we have
$$ g_{(n-p,p),(n-s,1^{s}),(n-t,1^t)}=\left\{\begin{array}{cl}
  2 & \ \mbox{if} \ \ t=s \ \mbox{and}\ s \geq p,\\1 & \ \mbox{if}\ \
  (t=s \ \mbox{and}\ s=p-1),\\  & \ \ \ \ \mbox{or}\
  (t=s+1 \ \mbox{and} \ s \geq p-1),  \\ & \ \ \ \ \mbox{or}\ (t=s-1\
  \mbox{and} \ s \geq p),\\ 0 & \ \mbox{otherwise.}
  \end{array} \right.$$
Similarly, Proposition~4.10 leads to formulas for the coefficient of a double hook
in $s_{(n-p,p)}*s_{(n-s, 1^s)}$.

\subsection{Multiplicities in $s_{(n-p,p)}*s_{\l}$}

In this section we use Theorem~3.2 to compute formulas for the multiplicity of Schur
functions corresponding to  two row partitions in the decomposition of
$s_{(n-p,p)}\ast s_\l$. Using this formula we show that under special conditions, if
$\l$ itself is a two row partition, these coefficients are unimodal. We have also
computed a formula for the coefficient of  $s_{(\n_1,\n_2,\n_3,\n_3)}$, $\n_3 \neq
0$,  in the decomposition of $s_{(n-p,p)}*s_{\l}$. The formulas are easy to program
and can yield the values of the coefficients for arbitrarily large $n$.

 \subsubsection{The multiplicity of $s_{(n-t,t)}$ in $s_{(n-p,p)}*s_{\l}$}
In this subsection we give a formula for the coefficient of $s_{(n-t,t)}$ in the product $s_{(n-p,p)}*s_{\l}$.
As corollaries of this formula we obtain simple formulas for the case when $\l$ is also a two row partition
and we show that under special conditions the Kronecker coefficients are unimodal.

\begin{prop}\label{th:coeff}
Let $n,t,p$ be non-negative integers, $p>0$, and $\l \vdash n$ be such that $\l_1
\geq 2p-1$. Let 
\begin{align*}
m_l&:=\min\{\l_1-\l_2, t-\l_2+p-2l-\max\{0,\l_3-l\}-\l_4, p-2l-1,\\
&\kern11cm 
\l_1+\l_4+p-2l-t\}\\
M_l&:=\max\{0, t-\l_2+p-2l-\l_3\}\\
m'_l&:=\min\{  \l_2-\max\{l,\l_3\},t- p+l-\max\{0,\l_3-l\}-\l_4,
\l_1-2p+3l,\\
&\kern10.3cm 
\l_4-t+\l_1-p+l+\l_2\}\\
M'_l&:= \max\{0,t-p+l-\l_3,\l_2-p+l\}\\
a_p&:= \chi(p \ \mbox{\em even})\chi\left(\l_3 \leq \frac{p}{2}\leq \min\{t,
\l_2\}\right)\chi(\l_2+\l_4\leq t \leq \min\{\l_2+\l_3, \l_1+\l_4\})\\
b(l)&:=\chi(\l_2-p+2l+\max\{0,\l_3-l\}+\l_4\leq t
\leq\l_2+\l_3-1)\\
c(l)&:= \chi(\l_1-\max\{p-l,\l_2\}\geq p-2l)
\chi(p-l+\max\{0,\l_2-p+l\}+\max\{0,\l_3-l\}\\
&\kern6.5cm
+\l_4\leq t \leq
\l_2+\l_3+p-2l-\max\{0,\l_3-l\})
\end{align*}

The   coefficient of $s_{(n-t,t)}$ in the decomposition of $s_{(n-p,p)}*s_{\l}$
equals 0 if $\ell(\l)>4$. If $\ell(\l)\leq 4$, then it equals 
\begin{multline*}
a_p+ \sum\limits_{l=\max\{\l_4,p-\l_2\}}^{\min\{\lfloor
\frac{p+1}{2}\rfloor -1, t, \l_2, p-\l_3\}} b(l) \max\{0,m_l-M_l+1\}\\+
\sum\limits_{l=\l_4}^{\min\{\lfloor \frac{p+1}{2} \rfloor-1,t, \l_2,p-\l_3\}} c(l)
\max\{0,m'_l-M'_l+1).
\end{multline*}
\end{prop}

\begin{proof}  We count Kronecker tableaux of shape $\l/\a$ and type $(n-t,t)/\a$, where
$\a \vdash p$, $\a \subseteq \l$ and $\a \subseteq (n-t,t)$. Thus $\a=(p-l,l)$ with
$l \leq \lfloor \frac{p}{2} \rfloor$, $l \leq \l_2$, $l \leq t$. Since the type of
the tableaux is $(n-t-p+l,t-l)$, the shape $\l/\a$ cannot contain three boxes
directly above each other. Thus, we  must have $\ell(\l) \leq 4$, $p-l \geq \l_3$
and $l \geq \l_4$. If $t<\l_4$, the coefficient of $s_{(n-t,t)}$ is zero. Hence, we
consider the Kronecker product $s_{(n-p,p)}*s_{(\l_1,\l_2,\l_3,\l_4)}$ with $\l_1
\geq 2p-1$, $\l_3\leq p$ and $\l_4 \leq t$. We have $k_{(p-l,l),(n-t,t)}^\l=0$
unless

\begin{equation}\label{lbound}\l_4\leq l\leq \min\{\lfloor
\frac{p}{2}\rfloor, t, \l_2, p-\l_3\}.\end{equation}

Moreover, whenever $\l/\a$ contains two boxes directly above each other, the top box
must be labelled $1$ and the bottom box must be labelled $2$.

\vskip10pt
\begin{picture}(400,50)
\put(100,0){\line(1,0){45}}\put(100,0){\line(0,0){40}}
\put(100,10){\line(1,0){115}}
\put(100,20){\line(1,0){160}}
\put(100,30){\line(1,0){185}}\put(100,40){\line(1,0){185}}
\put(102,2){\scriptsize{2}}\put(110,0){\line(0,0){20}}\put(102,12){\scriptsize{1}}
\put(115,2){$\cdots$}\put(115,12){$\cdots$}\put(135,0){\line(0,0){20}}
\put(137,2){\scriptsize{2}}\put(137,12){\scriptsize{1}}\put(145,0){\line(0,0){20}}
\put(150,12){$\cdots$}\put(170,10){\line(0,0){20}}\multiput(100,20)(.5,0){140}{\red{\line(0,0){10}}}
\put(172,12){\scriptsize{2}}\put(172,22){\scriptsize{1}}\put(180,10){\line(0,0){20}}
\put(185,12){$\cdots$}\put(185,22){$\cdots$}\put(205,10){\line(0,0){20}}
\put(207,12){\scriptsize{2}}\put(207,22){\scriptsize{1}}\put(215,10){\line(0,0){30}}
\multiput(100,30)(.5,0){230}{\red{\line(0,0){10}}}
\put(217,22){\scriptsize{2}}\put(217,32){\scriptsize{1}}
\put(225,20){\line(0,0){20}}\put(230,22){$\cdots$}\put(230,32){$\cdots$}
\put(250,20){\line(0,0){20}}\put(252,22){\scriptsize{2}}\put(252,32){\scriptsize{1}}
\put(260,20){\line(0,0){20}}\put(265,32){$\cdots$}\put(285,30){\line(0,0){10}}
\end{picture}
\vskip10pt

There must be enough $1$'s and $2$'s: 
\begin{align*}
n-t&\geq
p-l+\max\{0,\l_2-p+l\}+\max\{0,\l_3-l\}+\l_4\\
t&\geq
l+\max\{0,\l_2-p+l\}+\max\{0,\l_3-l\}+\l_4
\end{align*}

\noindent (I) Kronecker tableaux of shape $\l/(p/2,p/2)$ and type
$(n-t,t)/(p/2,p/2)$. Here $l=p/2$. Thus, $p$ must be even and since $p-l\geq \l_3$,
then $p/2\geq \l_3$. Condition (\ref{lbound}) becomes $\l_3\leq \frac{p}{2} \leq
\min\{t, \l_2\}$. The first row of $\l/(p-l,l)$ must be filled with $1$'s and the
second row must be filled with $2$'s. 

$\bullet$ Enough $1$'s: $n-t \geq
\l_1+\l_4$. 

$\bullet$ Enough $2$'s: $t \geq \l_2+\l_4$. 

$\bullet$ Lattice permutation: $t-\l_4 \leq \l_1$. 

\medskip
\noindent Hence,
\begin{multline*}
k_{(\frac{p}{2},\frac{p}{2}),(n-t,t)}^\l\\
 =\chi(p \ \mbox{even}) \chi\left(\l_3 \leq
\frac{p}{2}\leq \min\{t, \l_2\}\right)\chi(\l_2+\l_4\leq t \leq
\min\{\l_2+\l_3,\l_1+\l_4\}).
\end{multline*}
\medskip

\noindent (II) Kronecker tableaux of shape $\l/(p-l,l)$ and type $(n-t,t)/(p-l,l)$
with the box in position $(2, p-l)$ labelled $1$ and strictly less than $p-2l$ boxes
labelled $2$ in the first row. The tableaux with exactly $p-2l$ boxes labelled $2$
in the first row will be counted in (III). We consider only $l \neq p/2$, i.e., $l
\leq \lfloor \frac{p+1}{2}\rfloor -1$. 

\vskip10pt
\begin{picture}(400,50)
\put(100,0){\line(1,0){45}}\put(100,0){\line(0,0){40}}
\put(100,10){\line(1,0){115}}
\put(100,20){\line(1,0){205}}
\put(100,30){\line(1,0){285}}\put(100,40){\line(1,0){285}}
\put(102,2){\scriptsize{2}}\put(110,0){\line(0,0){20}}\put(102,12){\scriptsize{1}}
\put(115,2){$\cdots$}\put(115,12){$\cdots$}\put(135,0){\line(0,0){20}}
\put(137,2){\scriptsize{2}}\put(137,12){\scriptsize{1}}\put(145,0){\line(0,0){20}}
\put(150,12){$\cdots$}\put(170,10){\line(0,0){20}}\multiput(100,20)(.5,0){140}{\red{\line(0,0){10}}}
\put(172,12){\scriptsize{2}}\put(172,22){\scriptsize{1}}\put(180,10){\line(0,0){20}}
\put(185,12){$\cdots$}\put(185,22){$\cdots$}\put(205,10){\line(0,0){20}}
\put(207,12){\scriptsize{2}}\put(207,22){\scriptsize{1}}\put(215,10){\line(0,0){20}}
\multiput(100,30)(.5,0){320}{\red{\line(0,0){10}}}
\put(217,22){\scriptsize{1}}
\put(225,20){\line(0,0){10}}\put(230,22){$\cdots$}
\put(250,20){\line(0,0){10}}\put(252,22){\scriptsize{1}}
\put(260,20){\line(0,0){20}}
\put(270,20){\line(0,0){20}}\put(262,22){\scriptsize{2}}\put(262,32){\scriptsize{1}}
\put(275,22){$\cdots$}\put(275,32){$\cdots$}\put(295,20){\line(0,0){20}}
\put(305,20){\line(0,0){20}}\put(297,22){\scriptsize{2}}\put(297,32){\scriptsize{1}}
\put(310,32){$\cdots$}\put(330,30){\line(0,0){10}}\put(332,32){\scriptsize{1}}
\put(342,32){\scriptsize{2}}\put(340,30){\line(0,0){10}}\put(355,32){$\cdots$}
\put(350,30){\line(0,0){10}}\put(385,30){\line(0,0){10}}\put(377,32){\scriptsize{2}}
\put(375,30){\line(0,0){10}}
\end{picture}
\vskip10pt

$\bullet$ Position $(2, p-l)$ contained in $\l/(p-l,l)$: $l \geq p-\l_2$ (thus
$\max\{0,\l_2-p+l\}=\l_2-p+l$).

$\bullet$ Enough $1$'s (at most $p-2l-1$ boxes labelled $2$ in the first row):
$n-t\geq \l_1-p+2l+$

\ \ \ $1+p-2l+\l_4$, i.e., $t \leq \l_2+\l_3-1$.

$\bullet$ Enough $2$'s: $t\geq \l_2-p+2l+\max\{0,\l_3-l\}+\l_4$.
\medskip

We fill $k$ boxes in the first row with $2$'s.

$\bullet$ Strictly less than $p-2l$ boxes labelled $2$ in first row: $k \leq p-2l-1$

$\bullet$ Room for $2$'s: $k \leq \l_1-\l_2$.

$\bullet$ Total number of $2$'s: $k \leq t-\l_2+p-2l-\max\{0,\l_3-l\}-\l_4$

$\bullet$ Lattice permutation in first row: $k \leq p -2l$.

$\bullet$ Lattice permutation in third row: $ k \leq \l_1+\l_4+p-2l-t$.

$\bullet$ Minimum number of $2$'s in first row: $k \geq t-\l_2+p-2l-\l_3$
\medskip

Set
\begin{align*}
m_l&:=\min\{\l_1-\l_2, t-\l_2+p-2l-\max\{0,\l_3-l\}-\l_4, p-2l-1,\\
&\kern10cm
\l_1+\l_4+p-2l-t\}\\
M_l&:=\max\{0,t-\l_2+p-2l-\l_3\}.
\end{align*}

The number of Kronecker tableaux of shape $\l/(p-l,l)$ and type $(n-t,t)/(p-l,l)$
with the box in position $(2, p-l)$ labelled $1$ and strictly less than $p-2l$ boxes
labelled $2$ in the first row equals 
\begin{multline*}
 \sum_{l=\max\{\l_4,p-\l_2\}}^{\min\{\lfloor
\frac{p+1}{2} \rfloor-1, t, \l_2, p-\l_3\}} \chi(\l_2-p+2l+\max\{0,\l_3-l\}+\l_4\leq
t \leq\l_2+\l_3-1)\\
\cdot\max\{0, m_l-M_l+1\}.
\end{multline*}

\noindent (III) Kronecker tableaux of shape $\l/(p-l,l)$ and type $(n-t,t)/(p-l,l)$
with exactly $p-2l$ boxes in the first row labelled $2$. Again we consider only $l
\neq p/2$, i.e., $l \leq \lfloor \frac{p+1}{2}\rfloor -1$. 

\vskip10pt
\begin{picture}(400,50)
\put(100,0){\line(1,0){36}}\put(100,0){\line(0,0){40}}
\put(100,10){\line(1,0){100}}
\put(100,20){\line(1,0){208}}
\put(100,30){\line(1,0){285}}\put(100,40){\line(1,0){285}}
\put(102,2){\scriptsize{2}}\put(108,0){\line(0,0){20}}\put(102,12){\scriptsize{1}}
\put(112,2){$\cdots$}\put(112,12){$\cdots$}\put(128,0){\line(0,0){20}}
\put(130,2){\scriptsize{2}}\put(130,12){\scriptsize{1}}\put(136,0){\line(0,0){20}}
\put(140,12){$\cdots$}\put(156,10){\line(0,0){10}}\multiput(100,20)(.5,0){128}{\red{\line(0,0){10}}}
\put(164,10){\line(0,0){20}}\put(158,12){\scriptsize{1}}
\put(166,12){\scriptsize{2}}\put(166,22){\scriptsize{1}}\put(172,10){\line(0,0){20}}
\put(176,12){$\cdots$}\put(176,22){$\cdots$}\put(192,10){\line(0,0){20}}
\put(194,12){\scriptsize{2}}\put(194,22){\scriptsize{1}}\put(200,10){\line(0,0){20}}
\multiput(100,30)(.5,0){344}{\red{\line(0,0){10}}}
\put(202,22){\scriptsize{1}}\put(208,20){\line(0,0){10}}\put(212,22){$\cdots$}
\put(228,20){\line(0,0){10}}\put(230,22){\scriptsize{1}}\put(236,20){\line(0,0){10}}
\put(238,22){\scriptsize{2}}\put(244,20){\line(0,0){10}}\put(248,22){$\cdots$}
\put(264,20){\line(0,0){10}}\put(266,22){\scriptsize{2}}
\put(272,20){\line(0,0){20}}\put(280,20){\line(0,0){20}}
\put(274,22){\scriptsize{2}}\put(274,32){\scriptsize{1}}
\put(284,22){$\cdots$}\put(284,32){$\cdots$}\put(300,20){\line(0,0){20}}
\put(308,20){\line(0,0){20}}\put(302,22){\scriptsize{2}}\put(302,32){\scriptsize{1}}
\put(310,32){$\cdots$}\put(330,30){\line(0,0){10}}\put(332,32){\scriptsize{1}}
\put(342,32){\scriptsize{2}}\put(340,30){\line(0,0){10}}\put(355,32){$\cdots$}
\put(350,30){\line(0,0){10}}\put(385,30){\line(0,0){10}}\put(377,32){\scriptsize{2}}
\put(375,30){\line(0,0){10}}
\end{picture}
\vskip10pt

$\bullet$ Room for $2$'s in first row: $\l_1-\max\{p-l,\l_2\}\geq p-2l$.

$\bullet$ Enough $2$'s: $t\geq p-l+\max\{0,\l_2-p+l\}+\max\{0,\l_3-l\}+\l_4$.

$\bullet$ Enough $1$'s: $n-t \geq \l_1-p+2l+ \max\{0, \l_3-l\}+\l_4$.
\medskip

Fill $k$ boxes in the second row with  $2$'s.

$\bullet$ Room for $2$'s in second row and labels increasing in rows:

\ \ \ $\max\{0, \l_2-p+l\}\leq k \leq \l_2-\max\{l,\l_3\}$.

$\bullet$ Total number of $2$'s: $k \leq t- p+l-\max\{0,\l_3-l\}-\l_4$

$\bullet$ Lattice permutation in second row: $k \leq \l_1-2p+3l$.

$\bullet$ Lattice permutation in third row: $ k \leq \l_4-t+\l_1-p+l+\l_2$.

$\bullet$ Minimum number of $2$'s in second row: $k \geq t-p+l-\l_3$
\medskip

Set
\begin{align*}
m'_l&:=\min\{  \l_2-\max\{l,\l_3\},t- p+l-\max\{0,\l_3-l\}-\l_4,
\l_1-2p+3l,\\
&\kern10cm
\l_4-t+\l_1-p+l+\l_2\} \}\\
M'_l&:= \max\{0,t-p+l-\l_3,\l_2-p+l\}.
\end{align*}

The number of Kronecker tableaux of shape $\l/(p-l,l)$ and type $(n-t,t)/(p-l,l)$
with exactly $p-2l$ boxes in the first row labelled $2$ equals 
\begin{align*}
\sum_{l=\l_4}^{\min\{\lfloor \frac{p+1}{2} \rfloor-1, t, \l_2,
p-\l_3\}}&  \chi(\l_1-\max\{p-l,\l_2\}\geq p-2l)\\
&\cdot \chi(n-t \geq \l_1-p+2l+ \max\{0, \l_3-l\}+\l_4)\\ & \cdot\chi(t\geq
p-l+\max\{0,\l_2-p+l\}+\max\{0,\l_3-l\}+\l_4)\\ & \cdot\max\{0,m'_l-M'_l+1).
\end{align*}
\end{proof}

We now consider the special case when $\l$ is also a two row partition
in Theorem~\ref{th:coeff}. 
Let $\l=(n-s,s)\vdash n$ such that $n-s \geq 2p-1$.  In this case we
have $m_l=\min\{n-2s, t-s+p-2l, p-2l-1, n-s+p-2l-t\}$. Since $n-t\geq t$, we have
$m_l=\min\{n-2s, t-s+p-2l, p-2l-1\}$. Also $M_l=t-s+p-2l$. We see that $M_l \geq
m_l$. We obtain $m_l=M_l$ if and only if $t-s+p-2l\leq\min\{n-2s, p-2l-1\}$, i.e.,
$$m_l=M_l \ \mbox{if and only if} \ t \leq s-1 \ \mbox{and} \ l \geq 1/2(t+s+p-n).$$
Similarly, $m'_l=\min\{s-l, t-p+l, n-s-2p+3l\}$ and $M'_l = \max\{t-p+l, \max\{0,
s-p+l\}\}$. Again, $M'_l \geq m'_l$ and $$M'_l=m'_l \ \mbox{if and only if} \ l \geq
p-t, t \geq s, l \leq 1/2(p+s-t), l \geq 1/2(t+p+s-n).$$We also have
\begin{align*}
a_p&=\chi(p \mbox{\ even}) \chi\left(0 \leq \frac{p}{2}\leq
s\right)\chi(t=s)\\
b(l)&=\chi(t \leq s-1)\chi\left(l \leq \frac{t-s+p}{2}\right)\\
c(l)&=\chi(l \geq p-s)\chi(t \geq s) \chi\left(\frac{2s+p-n}{2}\leq l \leq
\frac{p+s-t}{2}\right)\\
&\kern2cm
+\chi(l<p-s)\chi\left(\max\left\{\frac{2p+s-n}{3},p-t\right\} \leq l \leq
\frac{p+s-t}{2}\right)
\end{align*}
We  obtain the following corollary.

\begin{cor}  Let $n,p,s,t$, be non-negative integers, such that $n-s \geq 2p-1$ and $p>0$. Let
\begin{alignat*}2
m_1&=\min\left\{ t, \left\lfloor\frac{t-s+p}{2}\right\rfloor\right\}& M_1&=\max\left\{0,p-s,
\left\lceil \frac{t+s+p-n}{2}\right\rceil\right\}\\ m_2&=\min\left\{s, \left\lfloor\frac{p+1}{2}\right\rfloor
-1\right\}& M_2&=\max\left\{0,p-s, \left\lceil \frac{2s+p-n}{2}\right\rceil\right\} \\
m_3&=\min\left\{ s, \left\lfloor\frac{p+s-t}{2}\right\rfloor\right\}& M_3&=M_1 \\ m_4&=\min\left\{ s,p-s-1,
\left\lfloor\frac{p+s-t}{2}\right\rfloor\right\}\kern.6cm&
M_4&=\max\left\{0,p-t, \left\lceil
\frac{t+s+p-n}{2}\right\rceil,\right.\\
&&&\kern4cm\left.
\left\lceil \frac{2p+s-n}{3}\right\rceil\right\}
\end{alignat*}

Then the coefficient of $s_{(n-t,t)}$ in the decomposition of
$s_{(n-p,p)}*s_{(n-s,s)}$ equals
$$
\begin{cases} \max\{0,m_1-M_1+1\} & \mbox{\ if\ } t<s,\\
\max\{0,m_2-M_2+1\}+ \chi(p \mbox{\ even})\chi(\frac{p}{2}\leq s)& \mbox{\ if\ } t=s,\\
\max\{0,m_3-M_3+1\} + \max\{0,m_4-M_4+1\} & \mbox{\ if\ } t>s.
\end{cases}
$$
\end{cor}

\begin{cor}  If $p \leq s$, $n-s \geq 2p-1$ and $n-p \geq 2s-1$. Then
 $$g_{(n-p,p),(n-s,s),(n-s,s)} =
  \begin{cases} \lfloor \frac{p}{2} \rfloor+1 & \mbox{if } n-p\geq 2s,
 \\
\lfloor \frac{p}{2} \rfloor & \mbox{if } n-p = 2s-1.\end{cases}
$$

\end{cor}

\begin{note} If $p \leq s$, $ 0 \leq l \leq \lfloor p/2 \rfloor$, there are
no Kronecker tableaux of shape $(n-s,s)/(p-l,l)$,  and type $(n-t,t)/(p-l,l)$ if
$t<s-p$ or $t>s+p$.
\end{note}

\begin{cor}
 If $p \leq s-1$, $n-s \geq 2p-1$, $n-p \geq 2s-1$ and $n-t \geq 2s-1$. Then
 the
coefficient of $s_{(n-t,t)}$  in $s_{(n-p,p)}\ast s_{(n-s,s)}$ is zero unless $s-p \leq t \leq s+p$.
If  $s-p \leq t \leq s+p$, then
$$g_{(n-p,p),(n-s,s),(n-t,t)}=\begin{cases} \lfloor \frac{p+t-s}{2} \rfloor +1
& \mbox{\ if \ } t\leq s-1,  \\   
 \lfloor \frac{p+s-t}{2}\rfloor +1 & \mbox{\ if \ } t\geq s \mbox{\ and \ } n-p \geq 2s,
 \\ \lfloor \frac{p+s-t}{2}\rfloor
 & \mbox{\ if \ } t\geq s  \mbox{\ and \ } n-p =2s-1, \end{cases}
$$ 
i.e., the
 sequence of coefficients, as $t=s-p,s-p+1, \ldots, s+p-1, s+p$ is unimodal. It is:
 $$\begin{array}{ll}
1, 1, 2, 2, \ldots, \lfloor \frac{p}{2}\rfloor +1, \ldots, 2,2,1,1& \mbox{\ if \ } n-p\geq
 2s,\\0, 0, 1, 1, \ldots, \lfloor \frac{p}{2}\rfloor , \ldots, 1,1,0,0& \mbox{\ if \ }
 n-p= 2s-1.\end{array}
$$
\end{cor}

\subsubsection{The multiplicity of $s_{(\n_1, \n_2, \n_3, \n_3)}$,   in $s_{(n-p,p)}*s_{(n-s,s)}$}

In this subsection we obtain a formula for $g_{(n-p,p),(n-s,s),(\n_1, \n_2, \n_3,
\n_3)}$, $\n_3 \neq 0$,  similar to the formulas given in \cite[Theorem~3.1 and
3.2]{rwh}. The advantage of our formula is that it does not involve cancellations.
 \begin{prop} Let $\n= (\n_1,\n_2,\n_3,\n_3)\vdash n$ with $\n_3
\neq 0$ and let $n,p,s$ be positive integers such that $n \geq 2p$, $n \geq 2s$ and
$n-s \geq 2p-1$. Let
\begin{align*}
M_1&:= \max\left\{p-\n_1,
\n_3,p-s+\n_3,\left\lceil\frac{1}{3}(2p-\n_1)\right\rceil,\left\lceil\frac{1}{2}(p+s-\n_3-\n_1)\right\rceil\right\}\\
m_1&:=\min\left\{\left\lfloor \frac{p+1}{2} \right\rfloor-1, \n_2,s-\n_3,
\left\lfloor\frac{1}{2}(\n_2+\n_3+p-s)\right\rfloor\right\}\\
M_2&:= \max\left\{\n_3, p-\n_2, \left\lceil
\frac{1}{2}(p+2s-n)\right\rceil\right\}\\
m_2&:=\min\left\{\left\lfloor \frac{p+1}{2}
\right\rfloor-1, \n_2,s-\n_3,\left\lfloor\frac{1}{2}(s+p-\n_2-\n_3)\right\rfloor\right\}
\end{align*}
The coefficient of $s_{\n}$ in the decomposition of $s_{(n-p,p)}*s_{(n-s,s)}$
equals
\begin{align*}
g_{(n-p,p),(n-s,s),\n} &=\chi(p \mbox{\
  even})\chi(\n_2+\n_3=s)\chi(\n_3 \leq p/2 \leq \n_2) \\ 
&\kern1cm  +\chi(\n_2+\n_3 \geq s) \max\{0,m_2-M_2+1\}  \\ & 
\kern1cm+\chi(\n_2+\n_3\leq s-1) \max\{0, m_1-M_1+1\} .
\end{align*}
\end{prop}

\begin{proof} Let $\a=(p-l,l)$ with $l \leq \lfloor p/2 \rfloor$, $l \leq s$, $l \leq
\n_2$, $l \geq p - \n_1$. We determine the number of Kronecker tableaux of shape
$(n-s,s)/(p-l,l)$ and type $\n/(p-l,l)$. In this special case, the number of $3$'s
equals the number of $4$'s. All $4$'s must be placed in the second row of $\l/\a$
and, because of the lattice permutation condition, all $3$'s must be placed in the
first row. We need: 
\medskip

\noindent $\bullet$ Room for $4$'s:  $s-l \geq \n_3$

\noindent $\bullet$ Room for $3$'s: $ n-s-p+l \geq \n_3$

\noindent $\bullet$ Lattice permutation in the first row: $l \geq \n_3$.

Since $n \geq 2p$ and $n \geq 2s$, the inequality  $l \geq \n_3$ implies $n-s-p+l
\geq \n_3$. Also, the inequality $\n_3 \leq l \leq s-\n_3$ implies $\n_3 \leq s/2$.

We consider only
 \begin{equation}\label{boundsl} \max\{p-\n_1,\n_3\} \leq l \leq
\min\{ \lfloor p/2\rfloor, s-\n_3, \n_2\}\end{equation}

 Note that the Kronecker coefficient of $s_{\n}$ is non-zero only if  $\n_3 \leq
\min\{p/2,s/2\}$.
\medskip

\noindent (I) Kronecker tableaux of shape $(n-p,p)/(p/2,p/2)$ and type
$\n/(p/2,p/2)$. In this case we must place all $1$'s  in the first row and all $2$'s
 in the second row of $(n-s,s)/(p-l,l)$. Therefore we need $\n_2+\n_3=s$.
 Condition (\ref{boundsl}) becomes
$$ \max\{p-\n_1,\n_3\} \leq p/2 \leq \min\{s-\n_3, \n_2 \}$$
Now, $p-\n_1 \leq p/2$ implies $\n_1 \geq p/2$. This condition is satisfied if $p/2
\leq \n_2$. Also, using $\n_2+\n_3 =s$, the inequality $p/2 \leq s-\n_3$ becomes
$p/2 \leq \n_2$. Thus the number of Kronecker tableaux of shape $(n-p,p)/(p/2,p/2)$
and type $\n/(p/2,p/2)$ equals
$$ \chi(p \mbox{\ even})\chi(\n_2+\n_3=s)\chi(\n_3 \leq p/2 \leq\n_2).$$

\noindent (II) Kronecker tableaux of shape $(n-p,p)/(p-l,l)$ and type $\n/(p-l,l)$,
$l \neq p/2$, such that the box in position $(2,p-l)$ is labelled $1$ and there are
strictly less than $p-2l$ boxes labelled $2$ in the first row. Since $l \neq p/2$,
we have $l \leq \lfloor \frac{p+1}{2}\rfloor -1$. In order to place label $1$ in
position $(2, p-l)$ we must have $p-l \leq s$. 

\vskip10pt
\begin{picture}(400,30)
\put(100,20){\line(1,0){255}}\put(100,0){\line(0,0){20}}
\put(100,10){\line(1,0){255}}\multiput(100,0)(.5,0){60}{\red{\line(0,0){10}}}
\put(100,0){\line(1,0){185}}\put(130,0){\line(0,0){10}}
\put(133,2){\scriptsize{1}}\put(140,0){\line(0,0){10}}
\put(145,2){$\cdots$}\put(165,0){\line(0,0){10}}
\put(168,2){\scriptsize{1}}\put(175,0){\line(0,0){20}}
\multiput(100,10)(.5,0){150}{\red{\line(0,0){10}}}
\put(178,2){\scriptsize{2}}\put(178,12){\scriptsize{1}}
\put(185,0){\line(0,0){20}}\put(190,2){$\cdots$}\put(190,12){$\cdots$}
\put(210,0){\line(0,0){20}}\put(213,2){\scriptsize{2}}\put(213,12){\scriptsize{1}}
\put(220,0){\line(0,0){20}}\put(223,2){\scriptsize{4}}\put(223,12){\scriptsize{1}}
\put(230,0){\line(0,10){20}}\put(235,12){$\cdots$}
\put(255,10){\line(0,0){10}}\put(258,12){\scriptsize{1}}\put(245,2){$\cdots$}
\put(265,10){\line(0,0){10}}\put(268,12){\scriptsize{2}}\put(275,10){\line(0,0){10}}
\put(280,12){$\cdots$}\put(300,10){\line(0,0){10}}\put(303,12){\scriptsize{2}}
\put(310,10){\line(0,0){10}}\put(313,12){\scriptsize{3}}\put(320,10){\line(0,0){10}}
\put(325,12){$\cdots$}\put(345,10){\line(0,0){10}}\put(347,12){\scriptsize{3}}
\put(355,10){\line(0,0){10}}\put(285,0){\line(0,0){10}}\put(275,0){\line(0,0){10}}
\put(278,2){\scriptsize{4}}
\end{picture}
\vskip10pt

There are $\n_3$ boxes labelled $4$ and $p-2l$ boxes labelled $1$ in the second row
and $\n_3$ boxes labelled $3$ in the first row. We need:

$\bullet$ Room for the $4$'s: $l \geq p-s+\n_3$.

$\bullet$ Enough $1$'s: $l \geq 1/3(2p-\n_1)$.

$\bullet$ Room for the $\n_1-2p+3l$ remaining  $1$'s in first row: $n-s-\n_3-p+l
\geq \n_1-2p+3l$,

\ \ \ i.e., $l \leq 1/2(n-s+p-\n_3-\n_1)=1/2(\n_2+\n_3+p-s)$.

$\bullet$ Number of $2$'s in the first row strictly less than $p-2l$:
$\n_2-(s-\n_3-p+2l)\leq p -2l-1$,

\ \ \ i.e., $\n_2+\n_3 \leq s-1$.

$\bullet$ Strictly increasing numbers in columns (no $2$'s above each other):
$s-\n_3 \leq \n_1-p+2l,$

\ \ \  i.e., $l \geq 1/2(p+s-\n_1-\n_3)$. 

Set 
\begin{align*}
M_1&:= \max\left\{p-\n_1,
\n_3,p-s+\n_3,\left\lceil\frac{1}{3}(2p-\n_1)\right\rceil,\left\lceil\frac{1}{2}(p+s-\n_3-\n_1)\right\rceil\right\},\\
m_1&:=\min\left\{\left\lfloor \frac{p+1}{2} \right\rfloor-1, \n_2,s-\n_3,
\left\lfloor\frac{1}{2}(\n_2+\n_3+p-s)\right\rfloor\right\}.
\end{align*}

The number of Kronecker tableaux of shape $(n-p,p)/(p-l,l)$ and type $\n/(p-l,l)$
such that the box in position $(2,p-l)$ is labelled $1$ and there are strictly less
than $p-2l$ boxes labelled $2$ in the first row equals
$$\chi(\n_2+\n_3\leq s-1) \max\{0, m_1-M_1+1\}.$$

\noindent (III) Kronecker tableaux of shape $(n-p,p)/(p-l,l)$ and type $\n/(p-l,l)$,
$l \neq p/2$,
  such
that exactly $p-2l$ boxes in the first row are labelled $2$. Again, $l \leq \lfloor
\frac{p+1}{2}\rfloor-1$. 

\vskip10pt
\begin{picture}(400,30)
\put(100,20){\line(1,0){255}}\put(100,0){\line(0,0){20}}
\put(100,10){\line(1,0){255}}\multiput(100,0)(.5,0){60}{\red{\line(0,0){10}}}
\put(100,0){\line(1,0){185}}\put(130,0){\line(0,0){10}}
\put(132,2){\scriptsize{1}}\put(138,0){\line(0,0){10}}
\put(142,2){$\cdots$}\put(158,0){\line(0,0){10}}
\put(160,2){\scriptsize{1}}\put(166,0){\line(0,0){10}}
\put(168,2){\scriptsize{2}}\put(174,0){\line(0,0){10}}
\put(195,2){$\cdots$}\put(202,10){\line(0,0){10}}
\multiput(100,10)(.5,0){204}{\red{\line(0,0){10}}}
\put(210,10){\line(0,0){10}}\put(204,12){\scriptsize{1}}
\put(214,12){$\cdots$}\put(230,0){\line(0,0){20}}
\put(238,0){\line(0,0){20}}\put(232,12){\scriptsize{1}}
\put(232,2){\scriptsize{2}}\put(240,2){\scriptsize{4}}
\put(246,0){\line(0,0){10}}\put(241,12){$\cdots$}
\put(256,10){\line(0,0){10}}\put(258,12){\scriptsize{1}}\put(254,2){$\cdots$}
\put(265,10){\line(0,0){10}}\put(268,12){\scriptsize{2}}\put(275,10){\line(0,0){10}}
\put(280,12){$\cdots$}\put(300,10){\line(0,0){10}}\put(303,12){\scriptsize{2}}
\put(310,10){\line(0,0){10}}\put(313,12){\scriptsize{3}}\put(320,10){\line(0,0){10}}
\put(325,12){$\cdots$}\put(345,10){\line(0,0){10}}\put(347,12){\scriptsize{3}}
\put(355,10){\line(0,0){10}}\put(285,0){\line(0,0){10}}\put(275,0){\line(0,0){10}}
\put(278,2){\scriptsize{4}}
\end{picture}
\vskip10pt

There are $\n_3$ boxes labelled $3$ and $p-2l$ boxes labelled $2$ in the first row
and $\n_3$ boxes labelled $4$ in the second row.  We need:

$\bullet$ Room for the $2$'s in first row: $n-s-\n_3-p+l \leq p-2l$, i.e., $l \geq
1/3(2p+s-n+\n_3)$.

\ \  Since $n-s \geq 2p-1$, we have $1/3(2p+s-n+\n_3) \leq 1/3(1+\n_3) \leq \n_3$
(here $\n_3 \neq 0$).

\ \  Since $l \geq \n_3$, the  condition $l \geq 1/3(2p+s-n+\n_3)$ is satisfied.

$\bullet$ Enough $2$'s: $l \geq p-\n_2$.

$\bullet$ Room for the $\n_2 -p +l$ remaining   $2$'s in the second row: $s-\n_3-l
\geq \n_2-p+l$,

\ \ \ i.e., $l\leq 1/2(s+p-\n_2-\n_3)$.

$\bullet$ Strictly increasing numbers in columns.

\ \ No $2$'s above each other: $s-\n_3 \leq n-s-p+2l-\n_3$, i.e., $l \geq
1/2(p+2s-n)$.

\ \ No $1$'s above each other: $s-\n_3-\n_2+p-l \leq p-l$, i.e., $\n_2+\n_3 \geq s$.
\medskip

Set 
\begin{align*}
M_2&:= \max\left\{\n_3, p-\n_2, \left\lceil \frac{1}{2}(p+2s-n)\right\rceil\right\},\\
m_2&:=\min\left\{\left\lfloor \frac{p+1}{2} \right\rfloor-1,
\n_2,s-\n_3,\left\lfloor\frac{1}{2}(s+p-\n_2-\n_3)\right\rfloor\right\}.
\end{align*}
The number of Kronecker tableaux of shape $(n-p,p)/(p-l,l)$ and type $\n/(p-l,l)$
such that exactly $p-2l$ boxes in the first row are labelled $2$ equals
$$\chi(\n_2+\n_3 \geq s) \max\{0, m_2-M_2+1\}.$$
\end{proof}

\begin{note} If $p=0$ or $p=1$ the multiplicity of $\n$ in Proposition~4.16
equals $0$.
\end{note}

\end{document}